\batchmode
\documentclass[11pt]{article}

\usepackage{epsfig}
\usepackage{graphicx}
\usepackage{color}
\usepackage{mathtools}

\newtheorem{theorem}{Theorem}
\newtheorem{lemma}{Lemma}
\newtheorem{corollary}{Corollary}

\newcommand{\be}{\begin{equation}}
\newcommand{\ee}{\end{equation}}
\newcommand{\bea}{\begin{eqnarray}}
\newcommand{\eea}{\end{eqnarray}}
\newcommand{\beas}{\begin{eqnarray*}}
\newcommand{\eeas}{\end{eqnarray*}}
\newcommand{\ba}{\begin{array}}
\newcommand{\ea}{\end{array}}

\definecolor{armygreen}{rgb}{0.29, 0.33, 0.13}

\newcommand{\real}{\mbox{$\mathbb{R}$}}

\newcommand{\Grad}{\ensuremath{\nabla}}
\newcommand{\eps}{\ensuremath{\epsilon}}

\newcommand{\bfa}{\ensuremath{\mathbf{a}}}
\newcommand{\bfb}{\ensuremath{\mathbf{b}}}

\newcommand{\bfe}{\ensuremath{\mathbf{e}}}

\newcommand{\bfi}{\ensuremath{\mathbf{i}}}
\newcommand{\bfj}{\ensuremath{\mathbf{j}}}

\newcommand{\bfq}{\ensuremath{\mathbf{q}}}
\newcommand{\bfr}{\ensuremath{\mathbf{r}}}

\newcommand{\bfI}{\ensuremath{\mathbf{I}}}
\newcommand{\bfB}{\ensuremath{\mathbf{B}}}

\newcommand{\bfF}{\ensuremath{\mathbf{F}}}

\def \bfvarphi{\mbox{\boldmath $\varphi$}}

\def\XXint#1#2#3{{\setbox0=\hbox{$#1{#2#3}{\int}$}
     \vcenter{\hbox{$#2#3$}}\kern-.5\wd0}}

\newcommand{\alpot}{\ensuremath{\frac{\alpha}{2}}}

\newcommand{\mcA}{\ensuremath{\mathcal{A}}}

\newcommand{\mcF}{\ensuremath{\mathcal{F}}}

\newcommand{\mcL}{\ensuremath{\mathcal{L}}}

\newcommand{\wtilde}{\ensuremath{\widetilde}}
\newcommand{\what}{\ensuremath{\widehat}}

\def\qed{\hbox{\vrule width 6pt height 6pt depth 0pt}}

\usepackage{amsmath}
\usepackage{amssymb}

\title{A generalized fractional Laplacian} 
\author{
	Xiangcheng Zheng\thanks{School of Mathematics, Shandong University, 
		Jinan 250100, China. email: {\tt xzheng@sdu.edu.cn}.} 
	\and
	V.J.~Ervin\thanks{School of Mathematical and Statistical Sciences,
	  Clemson University, Clemson, South Carolina 29634-0975, USA.
	  email: {\tt vjervin@clemson.edu}. }
	\and 
	 Hong Wang\thanks{Department of Mathematics, University of South Carolina, Columbia,
		South Carolina 29208, USA. email: {\tt hwang@math.sc.edu}. } 
}

\date{\today}

\begin{document}
\maketitle

\begin{abstract}
In this article we show that the fractional Laplacian in $\real^{2}$ can be factored into a product of the
divergence operator, a Riesz potential operator, and the gradient operator. Using this factored form we
introduce a generalization of the fractional Laplacian, involving a matrix $K(x)$, suitable when the 
fractional Laplacian is applied in a non homogeneous medium. For the case of $K(x)$ a constant, symmetric
positive definite matrix we show that the fractional Poisson equation is well posed, and determine
the regularity of the solution in terms of the regularity of the right hand side function.
\end{abstract}



\setcounter{equation}{0}
\setcounter{figure}{0}
\setcounter{table}{0}
\setcounter{theorem}{0}
\setcounter{lemma}{0}
\setcounter{corollary}{0}
\setcounter{definition}{0}
\section{Introduction}
 \label{sec_intro}
In recent years there has been considerable interest in studying nonlocal differential operators. The motivation
for this work has been to develop tools that can be used to more accurately model phenomena such as
anomalous diffusion, which cannot be suitable modeled using local operators. One such nonlocal operator
which has received much attention is the Fractional Laplace operator (a.k.a., the Fractional Laplacian). To
infer the Fractional Laplacian refers to a single operator is somewhat misleading as there are several
operators that come under this general heading. For example, for $0 < \alpha < 2$, on a bounded domain three 
such operators are \cite{aco171}: \\
(i) The Spectral Fractional Laplacian \\
The Spectral Fractional Laplacian, $\left( - \Delta \right)_{S}^{\alpot} \cdot$, is defined in terms of
the eigenvalues and eigenfunctions of the Laplace operator with homogeneous Dirichlet boundary
conditions on $\Omega$. Let $\left\{ \psi_{k} \right\}_{k \in \mathbb{N}}$ denote an orthonormal basis
for $L^{2}(\Omega)$ satisfying, for $\lambda_{k} \in \real^{+}$,
\begin{align*}
 - \Delta \psi_{k} &= \ \lambda_{k} \, \psi_{k} \ \mbox{ in } \Omega \, ,   \\
              \psi_{k} &= \ 0 \ \mbox{ on } \partial \Omega \, .
\end{align*}
Then $\left( - \Delta \right)_{S}^{\alpot} $ is defined for $u \in C_{0}^{\infty}(\Omega)$ as
\[ 
           \left( - \Delta \right)_{S}^{\alpot} u \ = \ \sum_{k = 1}^{\infty} \langle u \, , \, \psi_{k} \rangle \, \lambda_{k}^{\alpot} \, \psi_{k} \, ,
\]
 and by extension for $u \in H_{0}^{\alpot}(\Omega)$.

(ii) The Integral Fractional Laplacian \\
The Integral Fractional Laplacian, $\left( - \Delta \right)^{\alpot} \cdot$, is defined for $support(u) \subset \Omega \subset \real^{d}$
as
\be
 \left( - \Delta \right)^{\alpot} u(x) \ = \ \frac{1}{| \gamma_{d}(-\alpha)|} \, \lim_{\eps \rightarrow 0} 
 \int_{\real^{d} \backslash B(x, \eps)} \, \frac{u(x) \, - \, u(y)}{| x \, - \, y |^{d + \alpha}} \, dy \, , \ \ x \in \Omega \, ,
\label{defFLap}
\ee
where $\gamma_{d}(\alpha) \, := \, 2^{\alpha} \, \pi^{d/2} \, \Gamma(\alpot) / \Gamma( (d - \alpha)/2 )$, and 
$B(x, \eps)$ denotes the ball centered at $x$ with radius~$\eps$.

(iii) The Regional Fractional Laplacian \\
The Regional Fractional Laplacian is defined in a similar manner to \eqref{defFLap} but where the integration is restricted to 
$\Omega$ (i.e., $\real^{d} \backslash B(x, \eps)$ replaced by $\Omega \backslash B(x, \eps)$).

Another approach to defining the Fractional Laplacian, introduced by Caffarelli and Silvestre \cite{caf071}, is to consider
the Fractional Laplacian as the limit of a local operator problem defined on a semi-infinite cylinder $\Omega \times (0 , \infty)$.

Two recent papers which give insight into the relationship between various Fractional Laplace operators and their 
properties are \cite{kwa171, lis201}.

Related to the Integral Fractional Laplacian is the Riesz Potential operator, which is defined for $\alpha \in (0, d)$ by \cite{dyd171}
\[
   \left( - \Delta \right)^{- \alpot} u(x) \ = \ \frac{1}{\gamma_{d}(\alpha)} \, \int_{\real^{d}} 
   \frac{u(x - y)}{| y |^{d - \alpha}} \, dy \, .
\]

For suitably nice functions (see \cite[Proposition 1]{dyd171}) the Integral Fractional Laplacian is the inverse of the Riesz Potential 
operator.

In many physical applications ``motion'' occurs as a result of a local imbalance. For example, Fourier's Law and Fick's Law of Diffusion
state that the local flux is proportional to the gradient of the underlying quantity. For usual diffusion the impact of the
local imbalance is local. In the case of anomalous diffusion the impact of the local imbalance is nonlocal. For the
usual diffusion of the quantity $w$, the local flux in $\real^{d}$ is modeled as $K(x) \Grad w(x)$, and the spatial diffusion
operator as $\Grad \cdot K(x) \Grad w(x)$, where $K(x) \in \real^{d \times d}$ is a material parameter. (In the heat equation 
$K(x)$ represents the thermal conductivity of the medium.)

In order to include a material parameter into an anomalous diffusion model, consider (for $\bfq(x)$ denoting flux)
\be
    \begin{array}{ccc}
    \mbox{usual diffusion}    &  \hspace{1in}  &   \mbox{anomalous diffusion}   \\
   \Grad \cdot \bfq(x) \ = \ - \Grad \cdot K(x) \Grad w(x)  &      &   
     \Grad \cdot \bfq(x) \ = \  - \Grad \cdot \left(- \Delta \right)^{\frac{\alpha - 2}{2}} K(x) \Grad w(x) \, .    
    \end{array}
\label{defgenLap}
\ee

For $\bfI$ the identity matrix in $\real^{d \times d}$ and $k \in \real$, for $f$ a sufficiently nice function, a Fourier transform
argument shows that
\be
    k \left(- \Delta \right)^{\alpot} f \ = \ - \Grad \cdot \left(- \Delta \right)^{\frac{\alpha - 2}{2}} k \bfI \,  \Grad f \, .
\label{eqform1}
\ee 

The motivation for the work presented herein was to investigate a nonlocal fractional Laplace operator involving a material 
parameter $K(x)$. As a first step in this direction, we consider $K(x)$ to be a constant coefficient, symmetric, positive definite
matrix, and the domain $\Omega$ to be the unit disk in $\real^{2}$. More specifically, by the symmetry of the domain $\Omega$,
without loss of generality, we consider $K(x)$ to be a diagonal matrix with positive entries.
A corollary of our investigations (see Corollary \ref{cor4rep})
shows that \eqref{eqform1} holds for a family of functions having an edge singular behavior consistent with that 
expected of the solution of the Integral Fractional Laplacian. 

The problem investigated in this paper is: \textit{Given $f(x)$, $K(x) = \left[ \begin{array}{cc}
 k_{1}  &  0  \\  0  & k_{2}  \end{array} \right]$, $k_{1}, \, k_{2} \in \real^{+}$,
 determine $\wtilde{u}(x) \, = \, (1 - r^{2})_{+}^{\alpot} \, u(x)$ satisfying}
 \begin{align}
  \Grad \cdot \left(- \Delta \right)^{\frac{\alpha - 2}{2}} K(x) \Grad \wtilde{u}(x) &= \ f(x) \, , \ \ x \in \Omega \, ,   \label{deq1}  \\
       \wtilde{u}(x) &= \ 0  \, , \ \ x \in \real^{2} \backslash \Omega \, .   \label{deq2} 
 \end{align}
 
Note that \eqref{deq1} is different from the Riesz space fractional operator defined by (see \cite{zen141})
\[
  - k_{1} \, \frac{\partial^{\alpha} u}{\partial | x_{1} |^{\alpha}} \ - \ k_{2} \, \frac{\partial^{\alpha} u}{\partial | x_{2} |^{\alpha}} 
  \ = \ f(x_{1} , x_{2}) \, , \ \ (x_{1} , x_{2})  \in \Omega \, .
\]
This difference is apparent  when noting that $\left(- \Delta \right)^{\frac{\alpha - 2}{2}} \cdot $ is not a one dimensional
fractional integral operator.
For $\Omega \subset \real^{1}$ a detailed analysis of \eqref{deq1}--\eqref{deq2} is given in \cite{zhe211}.

Our analysis of \eqref{deq1}--\eqref{deq2} builds upon a result proven by Dyda et al. in \cite{dyd171}. For $L^{2}_{\beta}(\Omega)$
denoting the weighted $L^{2}$ space with domain $\Omega$, an orthogonal basis for $L^{2}_{\alpot}(\Omega)$ is
$\left\{ V_{l, \mu}(x) \, P_{n}^{(\alpot, l)}(2 r^{2} \, - \, 1) \right\}$, for $l, n = 0, 1, 2, \ldots, \infty$, $\mu \in \{-1, 1\}$, where
$V_{l, \mu}(x)$ denotes the solid harmonic polynomials in $\real^{2}$, and $P_{n}^{(a, b)}(t)$ the Jacobi polynomials
(see Section \ref{sec_prelim}). From \cite[Theorem 3]{dyd171},
\be
    \left(- \Delta \right)^{\alpot} (1 - r^{2})_{+}^{\alpot} \, V_{l, \mu}(x) \, P_{n}^{\alpot, l}(2 r^{2} \, - \, 1) 
    \ = \ \frac{2^{\alpha} \, \Gamma(n + 1 + \alpot) \, \Gamma(n + 1 + \alpot + l)}{\Gamma(n + 1) \, \Gamma(n + 1 + l)} \, 
     V_{l, \mu}(x) \, P_{n}^{\alpot, l}(2 r^{2} \, - \, 1) \, .
\label{dydres1}
\ee

A similar result to \eqref{dydres1} in the $\real^{1}$ setting is now well known and has been used to investigate the
existence, uniqueness and regularity of fractional order differential equations, as well as to develop spectral approximation
methods for such problems \cite{erv191, hao201, mao161, mao181, zhe211}.

In our analysis herein we show: 
\begin{enumerate}
\item A similar result to \eqref{dydres1} holds for the Riesz Potential operator, $\left(- \Delta \right)^{\frac{\alpha - 2}{2}} \cdot$
(see Theorem \ref{genThm3v2}).

\item For $u \in L^{2}_{\alpot}(\Omega)$, then 
$\Grad (1 - r^{2})^{\alpot} \, u \ \in \ (1 - r^{2})^{\alpot -1} \, \left[ \begin{array}{c}
  \text{span} \left\{ V_{l, \mu}(x) \, P_{n}^{\alpot - 1 \, ,  \, l}(2 r^{2} \, - \, 1) \right\}  \\
    \text{span} \left\{ V_{l, \mu}(x) \, P_{n}^{\alpot - 1 \, ,  \, l}(2 r^{2} \, - \, 1) \right\}  \end{array} \right]$
    (see Theorem \ref{gradinR2}). Notable is that this mapping holds
    when then the gradient is expressed in terms of the Cartesian coordinate system (i.e., in terms of $\bfi$ and $\bfj$), and
    not when written in the polar system (i.e., in terms of $\bfr$ and $\bfvarphi$).
    
\item (i) For  $\wtilde{u}(x) \, = \, (1 - r^{2})_{+}^{\alpot} \, u(x)$, and $u(x)$ expressed as a linear combination of the
basis functions $\left\{ V_{l, \mu}(x) \, P_{n}^{\alpot , l}(2 r^{2} \, - \, 1) \right\} $, substituting into \eqref{deq1} and
equating coefficients with the $f(x)$ (similarly expressed in terms of the basis), results in an infinite system of
coupled equations for the unknown coefficients (see \eqref{erq1}-\eqref{erq12}).

(ii) When the equations are appropriately ordered and rescaled, and the unknown coefficients renumbered and rescaled
the infinite system of equations decouple into tridiagonal block matrices of sizes $1, 1, 2, 2, \ldots, n, n, \ldots$ where the 
block matrices are uniformly invertible (see \eqref{Astr}-\eqref{Asub2}).

\item The Integral Fractional Laplacian can be factored into a product of the divergence operator, a Riesz Potential operator 
and the gradient operator (see Corollary \ref{cor4rep}). (Note that this result has been 
proven for the one-dimensional case (see, for example \cite{aco181, li211}), while, to the best of our knowledge, 
is not available in the literature for the two-dimensional operator.)

\item The dependence of the regularity of the solution of \eqref{deq1},\eqref{deq2} on the regularity of the RHS, $f(x)$.

\end{enumerate}

This paper is organized as follows. In the next section we introduce the Jacobi polynomials and the solid harmonic polynomials.
The combination of these polynomials form a basis for the weighted $L^{2}_{\beta}(\Omega)$ space, which plays a key role
in our analysis. The needed properties of $\Grad \cdot$ and $\left(- \Delta \right)^{\frac{\alpha - 2}{2}} \cdot$ applied
to functions $(1 - r^{2})^{\alpot} \, V_{l, \mu}(x) \, P_{n}^{(\alpot , l)}(2 r^{2} \, - \, 1)$ are derived in Section \ref{sec_MapR2}.
Using these results, in Section \ref{sec_ExUn} we establish the existence and uniqueness of solution,
$\wtilde{u}$, to \eqref{deq1}, \eqref{deq2}. The dependence of the regularity of $\wtilde{u}$ on the regularity of the
RHS, $f(x)$, is given in Section \ref{sec_Reg}. Contained in the Appendix are some ancillary results used in the analysis.

 \setcounter{equation}{0}
\setcounter{figure}{0}
\setcounter{table}{0}
\setcounter{theorem}{0}
\setcounter{lemma}{0}
\setcounter{corollary}{0}
\setcounter{definition}{0}
\section{Preliminaries}
 \label{sec_prelim}
In this section we present definitions and notation used in this article.
We follow the definitions and notation for the hypergeometric function used by
Dyda, Kuznetsov and Kwa\'{s}nicki in \cite{dyd171}.

\subsection{Hypergeometric Function}
\label{ssec_hypfun}
For $p, \, q$ nonnegative integers with $p \le \, q + 1$, and 
$\bfa \, = \, (a_{1}, \ldots , a_{p}) \in \mathbb{C}^{p}$, $\bfb \, = \, (b_{1}, \ldots , b_{q}) \in \mathbb{C}^{q}$,
with none of the $b_{j}$ nonpositive integers,
the generalized hypergeometric function is defined by
\be
_{p}F_{q}\left( \begin{array}{c}
                          \bfa \\  \bfb
                       \end{array} \, \Big{\vert} \, r \right) \ := \ 
             \sum_{n = 0}^{\infty} \frac{\overset{p}{\underset{j = 1}{\prod}} a_{j}^{(n)}}%
             {\overset{q}{\underset{j = 1}{\prod}} b_{j}^{(n)}} \, 
             \frac{r^{n}}{n !} \, ,
\label{defhypf}
\ee
where $c^{(n)} \ = \ c (c + 1) \cdots (c + n - 1)$ denoting the rising factorial, with $c^{(0)} \, = \, 1$.

Note that $_{p}F_{q}\left( \begin{array}{c}
                          \bfa \\  \bfb
                       \end{array} \, \Big{\vert} \, r \right)$ is invariant under the permutation of the
components of $\bfa$ or the components of $\bfb$.                

The regularized hypergeometric function, $_{p}\bfF_{q}(\cdot)$ is defined as
\be
_{p}\bfF_{q}\left( \begin{array}{c}
                          \bfa \\  \bfb
                       \end{array} \, \Big{\vert} \, r \right) \ := \ 
          \frac{1}{ \overset{q}{\underset{j = 1}{\prod}} \, \Gamma(b_{j}) } \, 
            _{p}F_{q}\left( \begin{array}{c}
                          \bfa \\  \bfb
                       \end{array} \, \Big{\vert} \, r \right)
\label{defghypf}
\ee

\subsection{Jacobi Polynomials}
\label{ssec_jpoly}
The Jacobi polynomials are defined as 
\begin{align}
P_{n}^{(a , b)}(t) &:= \ \frac{\Gamma(a \, + \, 1 \, + \, n)}{n !} \,
 _{2}\bfF_{1}\left( \begin{array}{c}
                          -n , \ 1 + a + b + n \\  a + 1  \end{array} \, \Big{\vert} \, \frac{1 - t}{2} \right)    \label{DKK40a} \\
&= \  \frac{(-1)^{n} \, \Gamma(b \, + \, 1 \, + \, n)}{n !} \,
 _{2}\bfF_{1}\left( \begin{array}{c}
                          -n , \ 1 + a + b + n \\  b + 1  \end{array} \, \Big{\vert} \, \frac{1 + t}{2} \right)    \label{DKK40b}     \\
&= \ (-1)^{n} \frac{\Gamma(n + 1 + b)}{n ! \, \Gamma(n + 1 + a + b)} \, \sum_{j = 0}^{n} (-1)^{j} \, 2^{-j} \, 
\left( \begin{array}{c}
           n  \\ j  \end{array} \right) \frac{\Gamma(n + j + 1 + a + b)}{\Gamma(j + 1 + b)} \, (1 + t)^{j} \, .  \label{Jsexp}                          
\end{align}

Note that $P_{0}^{(a , b)}(t) = 1$.

In case $a, \, b > -1$ the Jacobi polynomials satisfy the following orthogonality property.
\begin{align}
 & \int_{-1}^{1} (1 - t)^{a} (1 + t)^{b} \, P_{j}^{(a , b)}(t) \, P_{k}^{(a , b)}(t)  \, dt 
 \ = \
   \left\{ \begin{array}{ll} 
   0 , & k \ne j  \\
   |\| P_{j}^{(a , b)} |\|^{2}
   \, , & k = j  
    \end{array} \right.  \, ,  \nonumber \\
& \quad \quad \mbox{where } \  \ |\| P_{j}^{(a , b)} |\| \ = \
 \left( \frac{2^{(a + b + 1)}}{(2j \, + \, a \, + \, b \, + 1)} 
   \frac{\Gamma(j + a + 1) \, \Gamma(j + b + 1)}{\Gamma(j + 1) \, \Gamma(j + a + b + 1)}
   \right)^{1/2} \, .
  \label{spm22}
\end{align}  

We have the following differentiation properties
\begin{align}
\frac{d^{j}}{d t^{j}}  P_{k}^{(a , b)}(t)  &= \ 
\frac{\Gamma(a \, + \, b \, + \, k \, + \, 1 \, + \, j)}{2^{j} \, \Gamma(a \, + \, b \, + \, k \, + \, 1)} \, 
P_{k - j}^{(a+j \, , \, b+j)}(t)  \, ,   \label{diffJP}  \\
\frac{d}{d t} \left( (1 \, - \, t)^{\alpha/2} P_{k}^{(\alpot , m)}(t) \right) &= \ 
- (k \, + \, \alpot) (1 \, - \, t)^{\alpot - 1} \, P_{k}^{(\alpot -1  \, ,  \, m+1)}(t) \ 
\ \mbox{(see \cite[(E.1)]{hao211})} \, ,
\label{HLZZe1} \\
\frac{d}{d t} \left( (1 \, + \, t)^{m} P_{k}^{(\alpot , m)}(t) \right) &= \ 
 (k \, + \, m) (1 \, + \, t)^{m - 1} \, P_{k}^{(\alpot +1  \, ,  \, m-1)}(t) \, .
\label{HLZZe1b}
\end{align}

\subsection{Solid Harmonic Polynomials}
\label{ssec_shpoly}
The solid harmonic polynomials in $\real^{d}$ are the polynomials in $d$ variables which satisfy 
Laplace's equation.

In $\real^{2}$ the first few (linearly independent) solid harmonic polynomials are:
\[
1 , \ x_{1} , \  x_{2} , \  x_{1} x_{2} , \ x_{1}^{2} - x_{2}^{2}, \ldots 
\]

In $\real^{3}$ the first few (linearly independent) solid harmonic polynomials are:
\[
1 , \ x_{1} , \  x_{2} , \  x_{3} , \  x_{1} x_{2} , \  x_{1} x_{3} , \  x_{2} x_{3} , 
\ x_{1}^{2} - x_{2}^{2}, \ x_{1}^{2} - x_{3}^{2}, \ldots 
\]

 The solid harmonic polynomials of degree $l \ge 0$ form a finite dimensional vector space, having
 dimension
 \[
      M_{d , l} \ := \ \frac{d \, + \, 2 l \, - \, 2}{d \, + \, l \, - \, 2} 
      \left( \begin{array}{c}
             d \, + \, l \, - \, 2  \\   l  \end{array}  \right) \, .
\]

In $\real^{2}$ the solid harmonic polynomials of degree $l$ can be conveniently written in polar
coordinates, $\left( (r , \, \varphi ) \, : \, 0 \le r < \infty \, , \ 0 \le \varphi < \, 2 \pi \right)$, 
as $\{ r^{l} \cos (l \varphi ) \ , \   r^{l} \sin (l \varphi ) \} $. Note that $ M_{2 , l} = 2$.

\subsection{Function Spaces}
\label{ssec_funspc}     
The (weighted) $L^{2}(\Omega)$ space plays a central role in the analysis. For weight function $\beta(\cdot)$, $\beta(x) > 0 \, ,
x \in \Omega$, associated with $L^{2}_{\beta}(\Omega)$ we have the inner product and norm
\[
(f \, , g)_{\beta} \, := \, \int_{\Omega} \beta(x) \, f(x) \, g(x) \, d\Omega \, , \ \ \ \ 
\| f \|_{L^{2}_{\beta}(\Omega)} \, := \, (f \, , f)_{\beta}^{1/2} \, .
\]

In $\real^{d}$, $d = 2, 3$, $\Omega \, = \, \left\{ x \, : \, \| x \| < 1 \right\}$, i.e., the unit ball, a basis for 
$L^{2}_{\beta}(\Omega)$ is given as a product of the solid harmonic polynomials and Jacobi polynomials.

For $x = (r , \varphi) \in \Omega \subset \real^{2}$,  let
$ \omega^{\gamma} \, := \, (1 - r^{2})^{\gamma} $ .

\subsubsection{In $\real^{2}$}
\label{sssec_R2}     
\[
\mbox{Let } \ V_{l , 1}(x) \ := \ r^{l} \, \cos(l \varphi) \, , \  \ l = 0, 1, 2, \ldots  \ \ \mbox{ and } \
V_{l , -1}(x) \ := \ r^{l} \, \sin(l \varphi) \, , \  \ l = 1, 2, \ldots  .
\]

We also use the following notation
\[
    V_{l , \mu^{*}}(x) \ = \ \left\{ \begin{array}{rl}
    V_{l , -1}(x) & \mbox{  if } \mu = 1 \, ,   \\
    V_{l , 1}(x) & \mbox{  if } \mu = -1 \, .  \end{array} \right.
\]

Additionally, for a linear operator $\mcF \cdot$,  we use $\mcF(V_{l , \mu}(x)) \ = \ (\pm) V_{l , \sigma}(x)$ to denote
\[
\mcF(V_{l , \mu}(x)) = \ (\pm) V_{l , \sigma}(x) \ =  \ \left\{ \begin{array}{rl}
   + V_{l , \sigma}(x) & \mbox{  if } \mu = 1 \, , \\
   - V_{l , \sigma}(x) & \mbox{  if } \mu = -1 \, .  \end{array} \right.
\]
For example,
\[
 \frac{\partial}{\partial \varphi} V_{l , \mu}(x) \ = \ (\mp) V_{l , \mu^{*}}(x)
\]

In $\real^{2}$ an orthogonal basis for $L^{2}_{\beta}(\Omega)$ is \cite{dun141, li141}
\be
 \left\{  \cup_{l = 0}^{\infty} \cup_{n = 0}^{\infty} \left\{ V_{l , 1}(x) \, P_{n}^{(\beta  ,  l)}(2 r^{2} \, - \, 1) \right\} \right\}   \ \cup
  \left\{ \cup_{l = 1}^{\infty} \cup_{n = 0}^{\infty}  \left\{ V_{l , -1}(x) \, P_{n}^{(\beta  ,  l)}(2 r^{2} \, - \, 1) \right\} \right\}    \, .
\label{bsR2}
\ee

For notation brevity  we denote the basis in \eqref{bsR2} as
\[
  \cup_{l = 0}^{\infty} \cup_{n = 0}^{\infty}  \cup_{\mu = {1 , -1}} \left\{ V_{l , \mu}(x) \, P_{n}^{(\beta  ,  l)}(2 r^{2} \, - \, 1) \right\} 
\]
where we implicit assume that the terms   $V_{0, -1}(x) \, P_{n}^{(\beta  ,  l)}(2 r^{2} \, - \, 1), \, n = 0, 1, \ldots$ are omitted 
from the set.

\subsubsection{Action of the Riesz potential operator}
\label{sssec_RPO}     

The following theorem, used in our analysis below, extends a result in \cite{dyd171} from the fractional Laplace operator to the
Riesz potential operator.
\begin{theorem}  \label{genThm3v2}
For $\delta \ = \ d \, + \, 2 l$, $s$ an integer, $\alpot - s > -1$,
$f(x) \ = \ (1 \, - \, | x |^{2})^{\alpot - s}_{+} \, V_{l , \mu}(x) \, P_{n}^{(\alpot - s \, , \,  \frac{\delta}{2} - 1)}(2 |x|^{2} \, - \, 1)$,
\be
   \left( -\Delta \right)^{\frac{\alpha - 2}{2}} f(x) \ = \ (-1)^{1 - s} \,  2^{\alpha - 2} \, 
  \frac{ \Gamma(n + 1 - s + \frac{\alpha}{2}) \, \Gamma(n - 1 + \frac{\delta + \alpha}{2})}%
  {\Gamma(n + 1)\, \Gamma(n + 1 - s + \frac{\delta}{2})} \, 
  V_{l , \mu}(x) \, P_{n + 1 - s}^{(\alpot - 2 + s \, , \,  \frac{\delta}{2} - 1)}(2 |x|^{2} \, - \, 1) \, .
 \label{bceq1v}
\ee
\end{theorem}   
The proof of Theorem \ref{genThm3v2} begins by transforming $f(x)$, via a hypergeometric function, to a 
product of $V_{l , \mu}(x)$ and a Meijer G-function.
The operator $\left( -\Delta \right)^{\frac{\alpha - 2}{2}} \cdot$ is then applied and the
subsequent results transformed using properties of the Meijer G-function and the hypergeometric function
to \eqref{bceq1v}. As this is the only place in the paper where hypergeometric functions and Meijer G-functions
are used we present the proof of Theorem \ref{genThm3v2} in the Appendix.


 \setcounter{equation}{0}
\setcounter{figure}{0}
\setcounter{table}{0}
\setcounter{theorem}{0}
\setcounter{lemma}{0}
\setcounter{corollary}{0}
\setcounter{definition}{0}
\section{Mapping properties of $\Grad \cdot$ and $(- \Delta)^{\frac{\alpha - 2}{2}}$}
 \label{sec_MapR2}
In this section we investigate the mapping properties of  $\Grad \cdot$ and $(- \Delta)^{\frac{\alpha - 2}{2}}$
with the domain $\omega^{\alpot} \otimes L^{2}_{\alpot}(\Omega)$. Functions $g \in L^{2}_{\alpot}(\Omega)$ are
conveniently represented using polar coordinates but, surprisingly, the mapping properties of
 $\Grad \cdot$ and $(- \Delta)^{\frac{\alpha - 2}{2}}$ are far more revealing when written in a Cartesian framework
 (i.e., in terms of $\frac{\partial g}{\partial x}$ and $\frac{\partial g}{\partial y}$. To obtain this representation
 we firstly compute $\frac{\partial g}{\partial r}$ and $\frac{\partial g}{\partial \varphi}$, and then take
 an appropriate combination of these terms to obtain $\frac{\partial g}{\partial x}$ and $\frac{\partial g}{\partial y}$.
 
 \begin{theorem} \label{gradinR2}
 Let $f(x) \ = \ (1 - r^{2})^{\frac{\alpha}{2}} V_{l , \mu}(x) P_{n}^{(\frac{\alpha}{2} , l)}(2r^{2} - 1)$. Then for $n \ge 0$ and $l \ge 1$,
 \begin{align}
\frac{\partial f}{\partial x} &= \
 - (n + \alpot) (1 - r^{2})^{\alpot - 1} V_{l+1 , \mu}(x) P_{n}^{(\alpot-1 , l+1)}(2r^{2} - 1) \nonumber \\
 & \mbox{ } \quad \quad \ - \
(n + 1) (1 - r^{2})^{\alpot - 1} V_{l-1 , \mu}(x) P_{n+1}^{(\alpot-1 , l-1)}(2r^{2} - 1) \, ,    \label{ders1} \\
\frac{\partial f}{\partial y} &= \
 - ( \pm ) (n + \alpot) (1 - r^{2})^{\alpot - 1} V_{l+1 , \mu^{*}}(x) P_{n}^{(\alpot-1 , l+1)}(2r^{2} - 1) \nonumber \\
 & \mbox{ }  \quad  \quad   \ - \
( \mp ) (n + 1) (1 - r^{2})^{\alpot - 1} V_{l-1 , \mu^{*}}(x) P_{n+1}^{(\alpot-1 , l-1)}(2r^{2} - 1) \, .    \label{ders2}  
 \end{align}
 For the special case $l = 0$ and $\mu=1$, we have
 \begin{align}
\frac{\partial f}{\partial x} &= \
 -2 (n + \alpot) (1 - r^{2})^{\alpot - 1} V_{1 , 1}(x) P_{n}^{(\alpot-1 , 1)}(2r^{2} - 1) ,    \label{ders1z0} \\
\frac{\partial f}{\partial y} &= \
 - 2 (n + \alpot) (1 - r^{2})^{\alpot - 1} V_{1 , -1}(x) P_{n}^{(\alpot-1 , 1)}(2r^{2} - 1)  .    \label{ders2z0}  
 \end{align} 
 \end{theorem}
 
 The proof of Theorem \ref{gradinR2} is relatively long, using several of the recurrence formulas for Jacobi polynomials.
 So as not to unnecessarily distract from the existence and uniqueness of solution result for \eqref{deq1}, \eqref{deq2}, we give the
 proof of Theorem \ref{gradinR2} in the Appendix.


 \begin{theorem} \label{figradinR2}
 Let $f(x) \ = \ (1 - r^{2})^{\frac{\alpha}{2}} V_{l , \mu}(x) P_{n}^{(\frac{\alpha}{2} , l)}(2r^{2} - 1)$. Then for $n\geq 0$ and $l\geq 1$
 \begin{align}
(- \Delta)^{\frac{\alpha - 2}{2}} \, \frac{\partial f}{\partial x} &= \
  C_{2} \, (n + \alpot + l) \, V_{l+1 , \mu}(x) P_{n}^{(\alpot-1 , l+1)}(2r^{2} - 1) \nonumber \\
 & \mbox{ } \quad \quad \ + \
C_{2} \, (n + l + 1) \, V_{l-1 , \mu}(x) P_{n+1}^{(\alpot-1 , l-1)}(2r^{2} - 1) \, ,    \label{ders21} \\
(- \Delta)^{\frac{\alpha - 2}{2}} \, \frac{\partial f}{\partial y} &= \
  C_{2} \,  ( \pm ) (n + \alpot + l) \, V_{l+1 , \mu^{*}}(x) P_{n}^{(\alpot-1 , l+1)}(2r^{2} - 1) \nonumber \\
 & \mbox{ }  \quad  \quad   \ + \
 C_{2} \, ( \mp ) (n + l + 1) \, V_{l-1 , \mu^{*}}(x) P_{n+1}^{(\alpot-1 , l-1)}(2r^{2} - 1) \, ,    \label{ders22}  \\
\mbox{ where } C_{2} &= \  - 2^{\alpha - 2} \,   \frac{ \Gamma(n + 1 + \frac{\alpha}{2}) \, \Gamma(n + \frac{\alpha}{2} + l)}%
  {\Gamma(n + 1)\, \Gamma(n + 2 + l)}    \, .                      \label{ders23}
 \end{align}
 For the special case $l = 0$ and $\mu=1$, we have
  \begin{align}
(- \Delta)^{\frac{\alpha - 2}{2}} \, \frac{\partial f}{\partial x} &= \
  2C_{2} \, (n + \alpot ) \, V_{1 , 1}(x) P_{n}^{(\alpot-1 , 1)}(2r^{2} - 1) ,    \label{ders21z0} \\
(- \Delta)^{\frac{\alpha - 2}{2}} \, \frac{\partial f}{\partial y} &= \
  2C_{2} \,  (n + \alpot ) \, V_{1 ,-1}(x) P_{n}^{(\alpot-1 , 1)}(2r^{2} - 1)  .   \label{ders22z1}  
 \end{align}
 \end{theorem}
 
\textbf{Proof}: The stated results follow from applying Theorem \ref{genThm3v2} to the four relations in Theorem \ref{gradinR2}. 
\mbox{  } \hfill \qed

\begin{corollary} \label{gradwR2}
 Let $f(x) \ = \ V_{l , \mu}(x) P_{n}^{(\gamma , l)}(2r^{2} - 1)$. Then for $n \ge 1$ and $l \ge 1$,
 \begin{align}
\frac{\partial f}{\partial x} &= \
  (n + l)  \, V_{l-1 , \mu}(x) P_{n}^{(\gamma + 1 \, , \, l - 1)}(2r^{2} - 1) 
 \ + \
(n + \gamma + l + 1) \,  V_{l+1 , \mu}(x) P_{n-1}^{(\gamma + 1 \, , \, l + 1)}(2r^{2} - 1) \, ,    \label{dersw1} \\
\mbox{and} &  \nonumber \\
\frac{\partial f}{\partial y} &= \
 (\mp) (n + l) \,  V_{l-1 \, , \mu^{*}}(x) \, P_{n}^{(\gamma + 1 \, , \,  l - 1)}(2 r^{2} - 1)  
\ + \   (\pm) (n + \gamma + l + 1) \, V_{l+1 \, , \mu^{*}}(x) \,   P_{n-1}^{(\gamma + 1 \,  ,  \, l+1)}(2 r^{2} - 1)  \, .  
  \label{dersw2}  
 \end{align}
For $n = 0$ and $l \ge 1$,
\[
 \frac{\partial f}{\partial x}  \ = \ l  \, V_{l-1 , \mu}(x) P_{0}^{(\gamma  ,   l )}(2r^{2} - 1) \ \mbox{ and } \ \
 \frac{\partial f}{\partial y} \ = \ 
 (\mp) l \,  V_{l-1 \, , \mu^{*}}(x) \, P_{0}^{(\gamma  ,  l )}(2 r^{2} - 1)   \, .
\]
For $n \ge 1$ and $l = 0$, $\mu = 1$
\[
 \frac{\partial f}{\partial x}  \ = \ 2 (n + \gamma + 1)  \, V_{1 , 1}(x) P_{n-1}^{(\gamma + 1 \, , \, 1 )}(2r^{2} - 1) \ \mbox{ and } \ \
 \frac{\partial f}{\partial y} \ = \ 
 2 (n + \gamma + 1)  \, V_{1 , -1}(x) P_{n-1}^{(\gamma + 1 \, , \, 1 )}(2r^{2} - 1)  \, .
\]
For $n = 0$ and $l = 0$, $\mu = 1$
\[
 \frac{\partial f}{\partial x}  \ = \ 0 \ \mbox{ and } \ \
 \frac{\partial f}{\partial y} \ = \  0  \, .
\]
 \end{corollary}
\textbf{Proof}:  Using \eqref{diffJP} and $\frac{\partial f}{\partial x} \ = \ \cos(\varphi) \frac{\partial f}{\partial r} \, - \, 
\sin(\varphi) \, \frac{1}{r} \, \frac{\partial f}{\partial \varphi}$,
\begin{align*}
\frac{\partial f}{\partial x} &=  \ 
\cos(\varphi) \, \left( l \, V_{l , \mu}(x) \, r^{-1} \, P_{n}^{(\gamma ,  l)}(2 r^{2} - 1)  
\ + \ 
 2 (n + \gamma + l + 1) \, V_{l , \mu}(x) \, r \,  P_{n-1}^{(\gamma + 1 \,  ,  \, l + 1)}(2 r^{2} - 1) \right)  \\
& \quad  \quad - \ 
\sin(\varphi) \, \left(  \left( \frac{\partial}{\partial \varphi} V_{l , \mu}(x) \right) \, r^{-1} 
    \, P_{n}^{(\gamma ,  l)}(2 r^{2} - 1) \right) \, .   
\end{align*}
Then, using \eqref{ders5}, \eqref{ders6} and \eqref{ders7}, \eqref{ders8},
\begin{align*}
\frac{\partial f}{\partial x}  &= \
 l \, V_{l-1 , \mu}(x) \,  P_{n}^{(\gamma ,  l)}(2 r^{2} - 1)  \ + \ 
 (n + \gamma + l + 1) \, V_{l-1 , \mu}(x) \, r^{2} \,  P_{n-1}^{(\gamma + 1 \,  ,  \, l + 1)}(2 r^{2} - 1)  \\
& \quad \ + \ 
  (n + \gamma + l + 1) \, V_{l+1 , \mu}(x) \, P_{n-1}^{(\gamma + 1 \,  ,  \, l + 1)}(2 r^{2} - 1)  \\
&=  \
 (n + l) \, V_{l-1 , \mu}(x) \,  P_{n}^{(\gamma + 1 \, , \, l - 1 )}(2 r^{2} - 1)  \ + \ 
  (n + \gamma + l + 1) \, V_{l+1 , \mu}(x) \, P_{n-1}^{(\gamma + 1 \,  ,  \, l + 1)}(2 r^{2} - 1)   \, ,
\end{align*}
where, in the last step, we have used \eqref{vjJP1}.

Next, using $\frac{\partial f}{\partial y} \ = \ \sin(\varphi) \frac{\partial f}{\partial r} \, + \, 
\cos(\varphi) \, \frac{1}{r} \, \frac{\partial f}{\partial \varphi}$,
\begin{align*}
\frac{\partial f}{\partial y} &=  \ 
\sin(\varphi) \, \left( l \, V_{l , \mu}(x) \, r^{-1} \, P_{n}^{(\gamma ,  l)}(2 r^{2} - 1)  
\ + \ 
 2 (n + \gamma + l + 1) \, V_{l , \mu}(x) \, r \,  P_{n-1}^{(\gamma + 1 \,  ,  \, l + 1)}(2 r^{2} - 1) \right)  \\
& \quad  \quad + \ 
\cos(\varphi) \, \left(  \left( \frac{\partial}{\partial \varphi} V_{l , \mu}(x) \right) \, r^{-1} 
    \, P_{n}^{(\gamma ,  l)}(2 r^{2} - 1) \right) \, .   
\end{align*}
Using \eqref{ders15}, \eqref{ders16} and \eqref{ders17}, \eqref{ders18},
\begin{align*}
\frac{\partial f}{\partial y} &= \   
 l \, (\mp) V_{l-1 \, , \mu^{*}}(x) \, P_{n}^{(\gamma ,  l)}(2 r^{2} - 1)  
\ + \  (n + \gamma + l + 1) \, (\mp) V_{l-1 \, , \mu^{*}}(x) \, r^{2} \,  P_{n-1}^{(\gamma + 1 \,  ,  \, l+1)}(2 r^{2} - 1)   \\
& \quad \quad \ + \  (n + \gamma + l + 1) \, (\pm) V_{l+1 \, , \mu^{*}}(x) \,   P_{n-1}^{(\gamma + 1 \,  ,  \, l+1)}(2 r^{2} - 1)  \\
&= \  (\mp) (n + l) \,  V_{l-1 \, , \mu^{*}}(x) \, P_{n}^{(\gamma + 1 \, , \,  l - 1)}(2 r^{2} - 1)  
\ + \   (\pm) (n + \gamma + l + 1) \, V_{l+1 \, , \mu^{*}}(x) \,   P_{n-1}^{(\gamma + 1 \,  ,  \, l+1)}(2 r^{2} - 1) \, 
\end{align*}
where, again in the last step we have used \eqref{vjJP1}.  

The proofs of the other statements follow in a  similar manner and are therefore omitted for brevity. 
\mbox{ } \hfill \qed

\setcounter{equation}{0}
\setcounter{figure}{0}
\setcounter{table}{0}
\setcounter{theorem}{0}
\setcounter{lemma}{0}
\setcounter{corollary}{0}
\setcounter{definition}{0}
\section{Existence and Uniqueness of Solution}
 \label{sec_ExUn}
In this section we investigate the existence and uniqueness of solution to the following problem.
\textbf{Problem}: \\
Given $f$ and $k_{1}, k_{2} \in \real^{+}$, with $K(x) \, = \, \left[ \begin{array}{cc}
k_{1}  &   0 \\  0  &  k_{2}  \end{array} \right]$, find $\wtilde{u} \ = \ \omega^{\alpot} \, u$ such that
\begin{align}
\mcL(\wtilde{u})(x) \, := \, - \Grad \cdot (- \Delta)^{\frac{\alpha - 2}{2}} \, K(x) \, \Grad \wtilde{u}(x) &= \ f(x),  \ \ x \in \Omega\, , \label{plk1}  \\
\mbox{subject to   }  \wtilde{u} \ = \ 0 &\mbox{ on } \real^{2} \backslash \Omega \, .  \label{plk2}
\end{align}

We assume $u(x)$ can be expressed as 
\be
u(x) \, = \, \sum_{l\geq 1 , n\geq 0 , \mu\in \{1,-1\}} u_{l, n, \mu} \, V_{l , \mu}(x) \, P_{n}^{(\alpot , l)}(2 r^{2} - 1)+\sum_{n\geq 0}\frac{u_{0,n,1}}{2}V_{0 , 1}(x) \, P_{n}^{(\alpot , 0)}(2 r^{2} - 1)
\label{plk3}
\ee
and 
$$f(x) \, = \, \sum_{l\geq 1 , n\geq 0 , \mu\in \{1,-1\}} f_{l, n, \mu} \, V_{l , \mu}(x) \, P_{n}^{(\alpot , l)}(2 r^{2} - 1)+\sum_{n\geq 0} f_{0,n,1}V_{0 , 1}(x) \, P_{n}^{(\alpot , 0)}(2 r^{2} - 1)$$
 for coefficients
$u_{l, n, \mu}, \, f_{l, n, \mu} \in \real$. (Here we divide the coefficients $\{u_{0,n,1}\}_{n\geq 0}$ by $2$ in order that the matrices
$A_{3, 4, m}$ and $A_{3, 5, m}$ in \eqref{Asub1} and \eqref{Asub2} are symmetric.
Note that this factor of 2 results from the case of $l = 0, \, \mu = 1$ in Theorems \ref{gradinR2} and \ref{figradinR2}.)

Using Theorem \ref{figradinR2} and Corollary \ref{gradwR2}, $\mcL(\wtilde{u})$ can be expressed in terms of the basis
$\left\{ V_{l , \mu}(x) \, P_{n}^{(\alpot , l)}(2 r^{2} - 1) \right\}_{n, l, \mu}$. Equating the coefficients of $\mcL(\wtilde{u})$ with
those of $f$ we obtain a linear system of equations for the unknown coefficients $u_{l, n, \mu}$.

Using Theorem \ref{figradinR2},
$ k_{1} ( - \Delta )^{\frac{\alpha - 2}{2}} \frac{\partial \wtilde{u}}{\partial x} \ = \ T_{1} \ + \ T_{2}$, where 
\begin{align*}
 T_{1} &=  \ - 2^{\alpha - 2} \, k_{1} 
 \sum_{n \ge 0, \, \left\{   \scriptsize{\begin{array}{ll}
 l \ge 0, & \!  \!  \! \!  \!  \!  \mu = 1  \\
 l \ge 1,  & \!  \!  \!  \!  \!  \!  \mu = -1 \end{array} } \right.} u_{l, n, \mu} \, \frac{\Gamma(n + 1 + \alpot) \, \Gamma(n + 1 + \alpot + l)}%
{\Gamma(n + 1) \, \Gamma(n + 2 + l)} \, V_{l+1 \, , \, \mu}(x) \, P_{n}^{(\alpot - 1 \, , \, l+1)}(2 r^{2} - 1)   \\
&= \  - 2^{\alpha - 2} \, k_{1} 
 \sum_{n \ge 0, \, \left\{   \scriptsize{\begin{array}{ll}
 l \ge 1, & \!  \!  \! \!  \!  \!  \mu = 1  \\
 l \ge 2,  & \!  \!  \!  \!  \!  \!  \mu = -1 \end{array} } \right.} u_{l-1, n, \mu} \, \frac{\Gamma(n + 1 + \alpot) \, \Gamma(n + \alpot + l)}%
{\Gamma(n + 1) \, \Gamma(n + 1 + l)} \, V_{l , \mu}(x) \, P_{n}^{(\alpot - 1 \, , \, l)}(2 r^{2} - 1)  \, , \\
\end{align*}
and
\begin{align*}
 T_{2} &=  \ - 2^{\alpha - 2} \, k_{1} 
 \sum_{n \ge 0, \, \left\{   \scriptsize{\begin{array}{ll}
 l \ge 1, & \!  \!  \! \!  \!  \!  \mu = 1  \\
 l \ge 2,  & \!  \!  \!  \!  \!  \!  \mu = -1 \end{array} } \right.}
 u_{l, n, \mu} \, \frac{\Gamma(n + 1 + \alpot) \, \Gamma(n + \alpot + l)}%
{\Gamma(n + 1) \, \Gamma(n + 1 + l)} \, V_{l-1 \, , \, \mu}(x) \, P_{n+1}^{(\alpot - 1 \, , \, l-1)}(2 r^{2} - 1)   \\
&=  \ - 2^{\alpha - 2} \, k_{1} 
 \sum_{n \ge 1, \, \left\{   \scriptsize{\begin{array}{ll}
 l \ge 0, & \!  \!  \! \!  \!  \!  \mu = 1  \\
 l \ge 1,  & \!  \!  \!  \!  \!  \!  \mu = -1 \end{array} } \right.}
 u_{l+1, n-1, \mu} \, \frac{\Gamma(n + \alpot) \, \Gamma(n + \alpot + l)}%
{\Gamma(n) \, \Gamma(n + 1 + l)} \, V_{l , \, \mu}(x) \, P_{n}^{(\alpot - 1 \, , \, l)}(2 r^{2} - 1)%
\end{align*}

Using Corollary \ref{gradwR2},
\begin{align*}
\frac{\partial}{\partial x}T_{1} &= \ 
 - 2^{\alpha - 2} \, k_{1} 
\sum_{n = 0, \, \left\{   \scriptsize{\begin{array}{ll}
 l \ge 1, & \!  \!  \! \!  \!  \!  \mu = 1  \\
 l \ge 2,  & \!  \!  \!  \!  \!  \!  \mu = -1 \end{array} } \right.}
  u_{l-1, 0, \mu} \, \frac{\Gamma(0 + 1 + \alpot) \, \Gamma(0 + \alpot + l)}%
{\Gamma(0 + 1) \, \Gamma(0 + 1 + l)} \, l \,  V_{l-1 \, , \, \mu}(x) \, P_{0}^{(\alpot - 1 \, , \, l)}(2 r^{2} - 1)   \\
& \quad    - \, 2^{\alpha - 2} \, k_{1} 
\sum_{n \ge 1 , \, \left\{   \scriptsize{\begin{array}{ll}
 l \ge 1, & \!  \!  \! \!  \!  \!  \mu = 1  \\
 l \ge 2,  & \!  \!  \!  \!  \!  \!  \mu = -1 \end{array} } \right.}
  u_{l-1, n, \mu} \, \frac{\Gamma(n + 1 + \alpot) \, \Gamma(n + \alpot + l)}%
{\Gamma(n + 1) \, \Gamma(n + 1 + l)} \, \left( (n + l) \, V_{l-1 \,  , \, \mu}(x) \, P_{n}^{(\alpot  , \, l-1)}(2 r^{2} - 1)  \right. \\
& \quad \hspace{2in} 
+ \ \left. (n + \alpot + l) V_{l+1 \,  , \, \mu}(x) \, P_{n-1}^{(\alpot  , \, l+1)}(2 r^{2} - 1)  \right) \, .
\end{align*}

Recalling that $P_{0}^{(a , b)}(t) = 1$, we obtain
\begin{align}
\frac{\partial}{\partial x}T_{1} &= \ 
 - 2^{\alpha - 2} \, k_{1} 
\sum_{n = 0, \, \left\{   \scriptsize{\begin{array}{ll}
 l \ge 0, & \!  \!  \! \!  \!  \!  \mu = 1  \\
 l \ge 1,  & \!  \!  \!  \!  \!  \!  \mu = -1 \end{array} } \right.} 
u_{l, 0, \mu} \, \frac{\Gamma(1 + \alpot) \, \Gamma(1 + \alpot + l)}%
{\Gamma(1) \, \Gamma(1 + l)} \, V_{l ,  \mu}(x) \, 
  P_{0}^{(\alpot  , l)}(2 r^{2} - 1)   \nonumber  \\
& \ - \    2^{\alpha - 2} \, k_{1} 
\sum_{n \ge 1, \, \left\{   \scriptsize{\begin{array}{ll}
 l \ge 0, & \!  \!  \! \!  \!  \!  \mu = 1  \\
 l \ge 1,  & \!  \!  \!  \!  \!  \!  \mu = -1 \end{array} } \right.} 
u_{l, n, \mu} \, \frac{\Gamma(n + 1 + \alpot) \, \Gamma(n + 1 + \alpot + l)}%
{\Gamma(n + 1) \, \Gamma(n + 1 + l)} \,   V_{l  , \mu}(x) \, P_{n}^{(\alpot  , l)}(2 r^{2} - 1)  \nonumber   \\
& \ - \    2^{\alpha - 2} \, k_{1} 
\sum_{n \ge 0, \, \left\{   \scriptsize{\begin{array}{ll}
 l \ge 2, & \!  \!  \! \!  \!  \!  \mu = 1  \\
 l \ge 3,  & \!  \!  \!  \!  \!  \!  \mu = -1 \end{array} } \right.}  
u_{l-2, n+1, \mu} \, \frac{\Gamma(n + 2 + \alpot) \, \Gamma(n + 1 + \alpot + l)}%
{\Gamma(n + 2) \, \Gamma(n + 1 + l)}  \, V_{l  , \mu}(x) \, P_{n}^{(\alpot  , l)}(2 r^{2} - 1) \, .  \label{rew3} 
\end{align}
   
Next,
\begin{align}
\frac{\partial}{\partial x}T_{2} &= \ 
 - 2^{\alpha - 2} \, k_{1} 
\sum_{n \ge 1, \, l  = 0 ,  \, \mu = 1} u_{1, n-1, 1} \, \frac{\Gamma(n + \alpot) \, \Gamma(n + \alpot)}%
{\Gamma(n) \, \Gamma(n + 1)} \, 2 (n + \alpot) \,  V_{1 , 1}(x) \, P_{n-1}^{(\alpot  , 1)}(2 r^{2} - 1)  \nonumber  \\
& \quad    - \, 2^{\alpha - 2} \, k_{1} 
 \sum_{n \ge 1, \, \left\{   \scriptsize{\begin{array}{ll}
 l \ge 1, & \!  \!  \! \!  \!  \!  \mu = 1  \\
 l \ge 1,  & \!  \!  \!  \!  \!  \!  \mu = -1 \end{array} } \right.} 
u_{l+1, n-1, \mu} \, \frac{\Gamma(n  + \alpot) \, \Gamma(n + \alpot + l)}%
{\Gamma(n) \, \Gamma(n + 1 + l)} \, \left( (n + l) \, V_{l-1 \,  , \, \mu}(x) \, P_{n}^{(\alpot  , \, l-1)}(2 r^{2} - 1)  \right. \nonumber  \\
& \quad \hspace{2in} 
+ \ \left. (n + \alpot + l ) V_{l+1 \,  , \, \mu}(x) \, P_{n-1}^{(\alpot  , \, l+1)}(2 r^{2} - 1)  \right) \, \nonumber  \\
&= \  - 2^{\alpha - 2} \, k_{1} 
\sum_{n \ge 0, \, l  = 1 ,  \, \mu = 1} u_{1, n, 1} \, \frac{2 \, \Gamma(n + 1 + \alpot) \, \Gamma(n + 2 + \alpot)}%
{\Gamma(n + 1) \, \Gamma(n + 2)} \,   V_{1 , 1}(x) \, P_{n}^{(\alpot  , 1)}(2 r^{2} - 1)  \nonumber  \\
&  \quad   - 2^{\alpha - 2} \, k_{1} 
 \sum_{n \ge 1, \, \left\{   \scriptsize{\begin{array}{ll}
 l \ge 0, & \!  \!  \! \!  \!  \!  \mu = 1  \\
 l \ge 1,  & \!  \!  \!  \!  \!  \!  \mu = -1 \end{array} } \right.} 
u_{l+2, n-1, \mu} \, \frac{\Gamma(n + \alpot) \, \Gamma(n + 1 + \alpot + l)}%
{\Gamma(n) \, \Gamma(n + 1 + l)} \,   V_{l , \mu}(x) \, P_{n}^{(\alpot  , l)}(2 r^{2} - 1)  \nonumber  \\
&  \quad   - 2^{\alpha - 2} \, k_{1} 
\sum_{n \ge 0, \, \left\{   \scriptsize{\begin{array}{ll}
 l \ge 2, & \!  \!  \! \!  \!  \!  \mu = 1  \\
 l \ge 2,  & \!  \!  \!  \!  \!  \!  \mu = -1 \end{array} } \right.} 
u_{l, n, \mu} \, \frac{\Gamma(n + 1 + \alpot) \, \Gamma(n + 1 + \alpot + l)}%
{\Gamma(n + 1) \, \Gamma(n + 1 + l)} \,   V_{l , \mu}(x) \, P_{n}^{(\alpot  , l)}(2 r^{2} - 1)    \label{rew6} \, .
\end{align}

Again, using Theorem \ref{figradinR2},
$ k_{2} ( - \Delta )^{\frac{\alpha - 2}{2}} \frac{\partial \wtilde{u}}{\partial y} \ = \ T_{3} \ + \ T_{4}$, where 
\begin{align*}
 T_{3} &=  \ - 2^{\alpha - 2} \, k_{2} 
  \sum_{n \ge 0, \, \left\{   \scriptsize{\begin{array}{ll}
 l \ge 0, & \!  \!  \! \!  \!  \!  \mu = 1  \\
 l \ge 1,  & \!  \!  \!  \!  \!  \!  \mu = -1 \end{array} } \right.} 
 u_{l, n, \mu} \, \frac{\Gamma(n + 1 + \alpot) \, \Gamma(n + 1 + \alpot + l)}%
{\Gamma(n + 1) \, \Gamma(n + 2 + l)} \, (\pm) V_{l+1 \, , \, \mu^{*}}(x) \, P_{n}^{(\alpot - 1 \, , \, l+1)}(2 r^{2} - 1)   \\
&= \  - 2^{\alpha - 2} \, k_{2} 
 \sum_{n \ge 0, \, \left\{   \scriptsize{\begin{array}{ll}
 l \ge 1, & \!  \!  \! \!  \!  \!  \mu = 1  \\
 l \ge 2,  & \!  \!  \!  \!  \!  \!  \mu = -1 \end{array} } \right.} 
u_{l-1, n, \mu} \, \frac{\Gamma(n + 1 + \alpot) \, \Gamma(n + \alpot + l)}%
{\Gamma(n + 1) \, \Gamma(n + 1 + l)} \, (\pm) V_{l , \mu^{*}}(x) \, P_{n}^{(\alpot - 1 \, , \, l)}(2 r^{2} - 1)  \, , \\
\end{align*}
and
\begin{align*}
 T_{4} &=  \ - 2^{\alpha - 2} \, k_{2} 
  \sum_{n \ge 0, \, \left\{   \scriptsize{\begin{array}{ll}
 l \ge 2, & \!  \!  \! \!  \!  \!  \mu = 1  \\
 l \ge 1,  & \!  \!  \!  \!  \!  \!  \mu = -1 \end{array} } \right.}
 u_{l, n, \mu} \, \frac{\Gamma(n + 1 + \alpot) \, \Gamma(n + \alpot + l)}%
{\Gamma(n + 1) \, \Gamma(n + 1 + l)} \, (\mp) V_{l-1 \, , \, \mu^{*}}(x) \, P_{n+1}^{(\alpot - 1 \, , \, l-1)}(2 r^{2} - 1)   \\
&=  \ - 2^{\alpha - 2} \, k_{2} 
 \sum_{n \ge 1, \, \left\{   \scriptsize{\begin{array}{ll}
 l \ge 1, & \!  \!  \! \!  \!  \!  \mu = 1  \\
 l \ge 0,  & \!  \!  \!  \!  \!  \!  \mu = -1 \end{array} } \right.}
 u_{l+1, n-1, \mu} \, \frac{\Gamma(n + \alpot) \, \Gamma(n + \alpot + l)}%
{\Gamma(n) \, \Gamma(n + 1 + l)} \, (\mp) V_{l , \, \mu^{*}}(x) \, P_{n}^{(\alpot - 1 \, , \, l)}(2 r^{2} - 1)%
\end{align*}

Using Corollary \ref{gradwR2},
\begin{align}
\frac{\partial}{\partial y}T_{3} &= \ 
 - 2^{\alpha - 2} \, k_{2} 
 \sum_{n = 0, \, \left\{   \scriptsize{\begin{array}{ll}
 l \ge 1, & \!  \!  \! \!  \!  \!  \mu = 1  \\
 l \ge 2,  & \!  \!  \!  \!  \!  \!  \mu = -1 \end{array} } \right.}
u_{l-1, 0, \mu} \, \frac{\Gamma(0 + 1 + \alpot) \, \Gamma(0 + \alpot + l)}%
{\Gamma(0 + 1) \, \Gamma(0 + 1 + l)} \, l \,  V_{l-1 \, , \, \mu}(x) \, P_{0}^{(\alpot - 1 \, , \, l)}(2 r^{2} - 1)   \nonumber \\
& \quad    - \, 2^{\alpha - 2} \, k_{2} 
 \sum_{n \ge 1, \, \left\{   \scriptsize{\begin{array}{ll}
 l \ge 1, & \!  \!  \! \!  \!  \!  \mu = 1  \\
 l \ge 2,  & \!  \!  \!  \!  \!  \!  \mu = -1 \end{array} } \right.} 
u_{l-1, n, \mu} \, \frac{\Gamma(n + 1 + \alpot) \, \Gamma(n + \alpot + l)}%
{\Gamma(n + 1) \, \Gamma(n + 1 + l)} \, \left( (n + l) \, V_{l-1 \,  , \, \mu}(x) \, P_{n}^{(\alpot  , \, l-1)}(2 r^{2} - 1)  \right.  \nonumber \\
& \quad \hspace{2in} 
-  \ \left. (n + \alpot + l) V_{l+1 \,  , \, \mu}(x) \, P_{n-1}^{(\alpot  , \, l+1)}(2 r^{2} - 1)  \right) \,   \nonumber \\
 &= \ 
 - 2^{\alpha - 2} \, k_{2} 
 \sum_{n = 0, \, \left\{   \scriptsize{\begin{array}{ll}
 l \ge 0, & \!  \!  \! \!  \!  \!  \mu = 1  \\
 l \ge 1,  & \!  \!  \!  \!  \!  \!  \mu = -1 \end{array} } \right.} 
u_{l, 0, \mu} \, \frac{\Gamma(1 + \alpot) \, \Gamma(1 + \alpot + l)}%
{\Gamma(1) \, \Gamma(1 + l)} \, V_{l ,  \mu}(x) \, 
   P_{0}^{(\alpot  , l)}(2 r^{2} - 1)   \nonumber  \\
& \ - \    2^{\alpha - 2} \, k_{2} 
 \sum_{n \ge 1, \, \left\{   \scriptsize{\begin{array}{ll}
 l \ge 0, & \!  \!  \! \!  \!  \!  \mu = 1  \\
 l \ge 1,  & \!  \!  \!  \!  \!  \!  \mu = -1 \end{array} } \right.}
u_{l, n, \mu} \, \frac{\Gamma(n + 1 + \alpot) \, \Gamma(n + 1 + \alpot + l)}%
{\Gamma(n + 1) \, \Gamma(n + 1 + l)} \,   V_{l  , \mu}(x) \, P_{n}^{(\alpot  , l)}(2 r^{2} - 1)   \nonumber   \\
& \ + \    2^{\alpha - 2} \, k_{2} 
 \sum_{n \ge 0, \, \left\{   \scriptsize{\begin{array}{ll}
 l \ge 2, & \!  \!  \! \!  \!  \!  \mu = 1  \\
 l \ge 3,  & \!  \!  \!  \!  \!  \!  \mu = -1 \end{array} } \right.}
u_{l-2, n+1, \mu} \, \frac{\Gamma(n + 2 + \alpot) \, \Gamma(n + 1 + \alpot + l)}%
{\Gamma(n + 2) \, \Gamma(n + 1 + l)}  \, V_{l  , \mu}(x) \, P_{n}^{(\alpot  , l)}(2 r^{2} - 1) \, .  \label{rew13} 
\end{align}

Next,
\begin{align}
\frac{\partial}{\partial y}T_{4} &= \ 
 - 2^{\alpha - 2} \, k_{2} 
\sum_{n \ge 1 , \, l  = 0 ,  \, \mu = -1} 
u_{1, n-1, -1} \, \frac{\Gamma(n + \alpot) \, \Gamma(n + \alpot)}%
{\Gamma(n) \, \Gamma(n + 1)} \, (+)2 (n + \alpot) \,  V_{1 , -1}(x) \, P_{n-1}^{(\alpot  , 1)}(2 r^{2} - 1)  \nonumber  \\
& \quad    - \, 2^{\alpha - 2} \, k_{2} 
 \sum_{n \ge 1, \, \left\{   \scriptsize{\begin{array}{ll}
 l \ge 1, & \!  \!  \! \!  \!  \!  \mu = 1  \\
 l \ge 1,  & \!  \!  \!  \!  \!  \!  \mu = -1 \end{array} } \right.}
u_{l+1, n-1, \mu} \, \frac{\Gamma(n  + \alpot) \, \Gamma(n + \alpot + l)}%
{\Gamma(n) \, \Gamma(n + 1 + l)} \, \left( (-)(n + l) \, V_{l-1 \,  , \, \mu}(x) \, P_{n}^{(\alpot  , \, l-1)}(2 r^{2} - 1)  \right. \nonumber  \\
& \quad \hspace{2in} 
+ \ \left. (n + \alpot + l ) V_{l+1 \,  , \, \mu}(x) \, P_{n-1}^{(\alpot  , \, l+1)}(2 r^{2} - 1)  \right) \, \nonumber  \\
&= \  - 2^{\alpha - 2} \, k_{2} 
\sum_{n \ge 0 , \, l  = 1 ,  \, \mu = -1} u_{1, n, -1} \, \frac{2 \, \Gamma(n + 1 + \alpot) \, \Gamma(n + 2 + \alpot)}%
{\Gamma(n + 1) \, \Gamma(n + 2)} \,   V_{1 , -1}(x) \, P_{n}^{(\alpot  , 1)}(2 r^{2} - 1)   \nonumber  \\
&  \quad   + \,  2^{\alpha - 2} \, k_{2} 
 \sum_{n \ge 1, \, \left\{   \scriptsize{\begin{array}{ll}
 l \ge 0, & \!  \!  \! \!  \!  \!  \mu = 1  \\
 l \ge 0,  & \!  \!  \!  \!  \!  \!  \mu = -1 \end{array} } \right.}
u_{l+2, n-1, \mu} \, \frac{\Gamma(n + \alpot) \, \Gamma(n + 1 + \alpot + l)}%
{\Gamma(n) \, \Gamma(n + 1 + l)} \,   V_{l , \mu}(x) \, P_{n}^{(\alpot  , l)}(2 r^{2} - 1)   \nonumber  \\
&  \quad   - \, 2^{\alpha - 2} \, k_{2} 
 \sum_{n \ge 0, \, \left\{   \scriptsize{\begin{array}{ll}
 l \ge 2, & \!  \!  \! \!  \!  \!  \mu = 1  \\
 l \ge 2,  & \!  \!  \!  \!  \!  \!  \mu = -1 \end{array} } \right.} 
u_{l, n, \mu} \, \frac{\Gamma(n + 1 + \alpot) \, \Gamma(n + 1 + \alpot + l)}%
{\Gamma(n + 1) \, \Gamma(n + 1 + l)} \,   V_{l , \mu}(x) \, P_{n}^{(\alpot  , l)}(2 r^{2} - 1)  \label{rew16} \, .
\end{align}


As mentioned at the beginning of this section, we obtain the determining equations for the unknown coefficients $u_{l, n, \mu}$ by equating
the coefficients of the basis functions $\left\{ V_{l , \mu}(x) \, P_{n}^{(\alpot , l)}(2 r^{2} - 1) \right\}_{n, l, \mu}$ for
$\mcL(\wtilde{u}) \, = \, - \Grad \cdot K(x) (-\Delta)^{\frac{\alpha - 2}{2}}) \Grad \wtilde{u}$ with those of $f$.

\begin{figure}[!ht]
\begin{minipage}{.46\linewidth}
\begin{center}
 \includegraphics[height=2.25in]{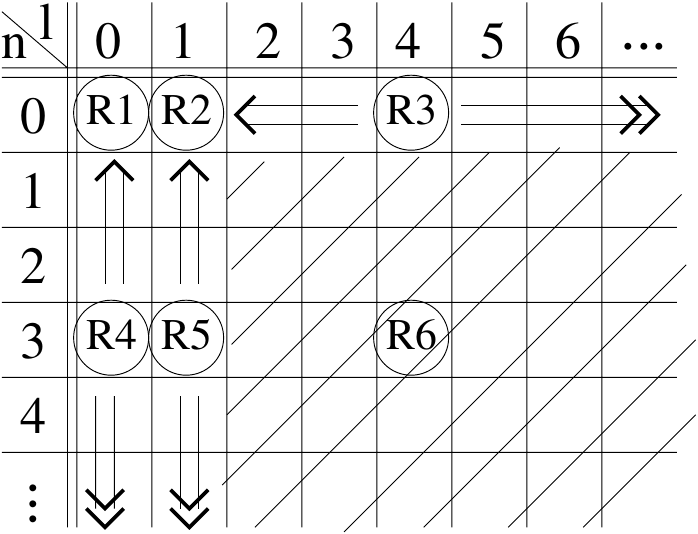}
 \caption{Partition of the index set $\{(n , l)\}_{n \ge 0 , l \ge 0}$ for $\mu = 1$ into $R1, R2, \ldots, R6$.}
\label{part1}
\end{center}

\end{minipage} \hfill
\begin{minipage}{.46\linewidth}
 
\begin{center}
 \includegraphics[height=2.25in]{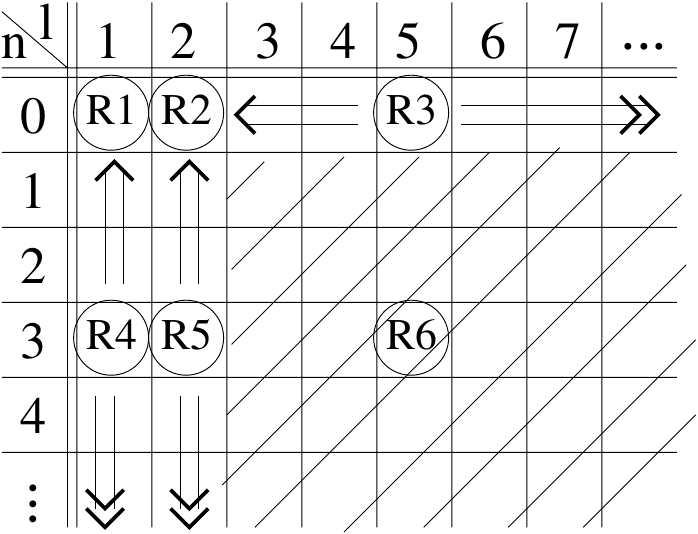}
 \caption{Partition of the index set $\{(n , l)\}_{n \ge 0 , l \ge 1}$ for $\mu = -1$ into $R1, R2, \ldots, R6$.}
\label{part2}
\end{center}
\end{minipage} 
\end{figure}

To obtain the determining equations it is convenient to group the associated basis functions with the index sets $\{ (n, l) \}_{n \ge 0, \, l \ge 0}$
for $\mu = 1$, and $\{ (n, l) \}_{n \ge 0, \, l \ge 1}$ for $\mu = -1$, separately. We then partition the index sets into 6 pieces as
follows (see Figures \ref{part1} and \ref{part2}):
\[
\begin{array}{clccl}
\mbox{For } \mu = 1 :  & R1 \, = \, \{ (0 , 0) \} \, ;  & \quad \quad  & \mbox{For } \mu = -1 :  & R1 \, = \, \{ (0 , 1) \} \, ;  \\
   & R2 \, = \, \{ (0 , 1) \} \, ;  & \quad \quad  &    & R2 \, = \, \{ (0 , 2) \} \, ;  \\
   & R3 \, = \, \{ (0 , l) \}_{l \ge 2} \, ;  & \quad \quad  &    & R3 \, = \, \{ (0 , l) \}_{l \ge 3} \, ;  \\
   & R4 \, = \, \{ (n , 0) \}_{n \ge 1} \, ;  & \quad \quad  &    & R4 \, = \, \{ (n , 1) \}_{n \ge 1} \, ;  \\
   & R5 \, = \, \{ (n , 1) \}_{n \ge 1} \, ;  & \quad \quad  &    & R5 \, = \, \{ (n , 2) \}_{n \ge 1} \, ;  \\
   & R6 \, = \, \{ (n , l) \}_{n \ge 1, \, l \ge 2} \, ;  & \quad \quad  &    & R6 \, = \, \{ (n , l) \}_{n \ge 1, \, l \ge 3} \, .
\end{array}
\]

Combining \eqref{rew3}-\eqref{rew16} and equating with the coefficients of 
$f(x) \, = \, \sum_{n \, l \, \mu} f_{l, n, \mu} \, V_{l , \mu}(x) \, P_{n}^{(\alpot , l)}(2 r^{2} - 1)$ we obtain the following determining system of
equations for the unknowns $u_{l, n, \mu}$.

\underline{For $\mu = 1$}: 
\begin{align}
 R1 \, (n = 0 , l = 0):  \quad & 2^{\alpha - 2} \, (k_{1} + k_{2}) \, \frac{\Gamma(1 + \alpot) \, \Gamma(1 + \alpot)}{\Gamma(1) \, \Gamma(1)} \,
  u_{0, 0, 1}  \ = \ f_{0, 0, 1}      \label{erq1}  \\
R2 \, (n = 0 , l = 1): \quad & 2^{\alpha - 2} \, \left( 2 k_{1} \, + \, (k_{1} + k_{2}) \right) \, 
\frac{\Gamma(1 + \alpot) \, \Gamma(2 + \alpot)}{\Gamma(1) \, \Gamma(2)} \,
  u_{1, 0, 1}  \ = \ f_{1, 0, 1}      \label{erq2}  \\
R3 \, \{ (n = 0 , l ) \}_{l \ge 2}: \quad &   \nonumber  \\
 & \hspace{-1.75in} 2^{\alpha - 2} \, \left( 2 (k_{1} + k_{2}) \, 
\frac{\Gamma(1 + \alpot) \, \Gamma(1 + \alpot + l)}{\Gamma(1) \, \Gamma(1 + l)} \,
  u_{l, 0, 1}   \ +  \ 
  (k_{1} - k_{2}) \, 
\frac{\Gamma(2 + \alpot) \, \Gamma(1 + \alpot + l)}{\Gamma(2) \, \Gamma(1 + l)} \,
  u_{l-2, 1, 1}
  \right) \nonumber \\
 & \hspace{3.5in}  \ = \ f_{l, 0, 1}      \label{erq3}  \\
R4 \, \{ (n , l = 0) \}_{n \ge 1}: \quad &   \nonumber  \\
 & \hspace{-1.75in} 2^{\alpha - 2} \, \left( (k_{1} - k_{2}) \, 
\frac{\Gamma(n + \alpot) \, \Gamma(n + 1 + \alpot)}{\Gamma(n) \, \Gamma(n + 1)} \,
  u_{2, n-1, 1}   \ +  \ 
  (k_{1} + k_{2}) \, 
\frac{\Gamma(n + 1 + \alpot) \, \Gamma(n + 1 + \alpot)}{\Gamma(n + 1) \, \Gamma(n + 1)} \,
  u_{0, n, 1}
  \right) \nonumber \\
  & \hspace{3.5in} \ = \ f_{0, n, 1}      \label{erq4}  \\
R5 \, \{ (n , l = 1) \}_{n \ge 1}: \quad &   \nonumber  \\
 & \hspace{-1.75in} 2^{\alpha - 2} \, \left( (k_{1} - k_{2}) \, 
\frac{\Gamma(n + \alpot) \, \Gamma(n + 2 + \alpot)}{\Gamma(n) \, \Gamma(n + 2)} \,
  u_{3, n-1, 1}   \ +  \ 
(2 k_{1} \, + \,   (k_{1} + k_{2}) ) \, 
\frac{\Gamma(n + 1 + \alpot) \, \Gamma(n + 2 + \alpot)}{\Gamma(n + 1) \, \Gamma(n + 2)} \,
  u_{1, n, 1}
  \right)  \nonumber \\
 & \hspace{3.5in}  \ = \ f_{1, n, 1}      \label{erq5}  \\
R6 \, \{ (n , l) \}_{n \ge 1 , l \ge 2}: \quad &   \nonumber  \\
 & \hspace{-1.75in} 2^{\alpha - 2} \, \left( (k_{1} - k_{2}) \, 
\frac{\Gamma(n + \alpot) \, \Gamma(n + 1 + \alpot + l)}{\Gamma(n) \, \Gamma(n + 1 + l)} \,
  u_{l+2, n-1, 1}   \ +  \ 
2 (k_{1} + k_{2})  \, 
\frac{\Gamma(n + 1 + \alpot) \, \Gamma(n + 1 + \alpot + l )}{\Gamma(n + 1) \, \Gamma(n + 1 + l)} \,
  u_{l, n, 1}   \right.    \nonumber  \\
 & \hspace{0.2in} \left. + \  (k_{1} - k_{2}) \, 
\frac{\Gamma(n + 2 + \alpot) \, \Gamma(n + 1 + \alpot + l)}{\Gamma(n + 2) \, \Gamma(n + 1 + l)} \,
  u_{l-2, n+1, 1}
  \right) \ = \ f_{l, n, 1}    \, .  \label{erq6}  
\end{align}

\underline{For $\mu = -1$}: 
\begin{align}
 R1 \, (n = 0 , l = 1):  \quad & 2^{\alpha - 2} \, \left( 2 k_{2} \, + \, (k_{1} + k_{2}) \right) \, 
 \frac{\Gamma(1 + \alpot) \, \Gamma(2 + \alpot)}{\Gamma(1) \, \Gamma(2)} \,
  u_{1, 0, -1}  \ = \ f_{1, 0, -1}      \label{erq7}  \\
R2 \, (n = 0 , l = 2): \quad & 2^{\alpha - 2} \  2 (k_{1} + k_{2}) \, 
\frac{\Gamma(1 + \alpot) \, \Gamma(3 + \alpot)}{\Gamma(1) \, \Gamma(3)} \,
  u_{2, 0, -1}  \ = \ f_{2, 0, -1}      \label{erq8}  \\
R3 \, \{ (n = 0 , l) \}_{l \ge 3}: \quad &   \nonumber  \\
 & \hspace{-1.75in} 2^{\alpha - 2} \, \left( 2 (k_{1} + k_{2}) \, 
\frac{\Gamma(1 + \alpot) \, \Gamma(1 + \alpot + l)}{\Gamma(1) \, \Gamma(1 + l)} \,
  u_{l, 0, -1}   \ +  \ 
  (k_{1} - k_{2}) \, 
\frac{\Gamma(2 + \alpot) \, \Gamma(1 + \alpot + l)}{\Gamma(2) \, \Gamma(1 + l)} \,
  u_{l-2, 1, -1}
  \right) \nonumber \\
 & \hspace{3.5in}  \ = \ f_{l, 0, -1}      \label{erq9}  \\
R4 \, \{ (n , l = 1) \}_{n \ge 1}: \quad &   \nonumber  \\
 & \hspace{-1.75in} 2^{\alpha - 2} \, \left( (k_{1} - k_{2}) \, 
\frac{\Gamma(n + \alpot) \, \Gamma(n + 2 + \alpot)}{\Gamma(n) \, \Gamma(n + 2)} \,
  u_{3, n-1, -1}   \ +  \ 
(2 k_{2} +  (k_{1} + k_{2}) ) \, 
\frac{\Gamma(n + 1 + \alpot) \, \Gamma(n + 2 + \alpot)}{\Gamma(n + 1) \, \Gamma(n + 2)} \,
  u_{1, n, 1}
  \right)   \nonumber \\
  & \hspace{3.5in} \ = \ f_{1, n, -1}      \label{erq10}  \\
R5 \, \{ (n , l = 2) \}_{n \ge 1}: \quad &   \nonumber  \\
 & \hspace{-1.75in} 2^{\alpha - 2} \, \left( (k_{1} - k_{2}) \, 
\frac{\Gamma(n + \alpot) \, \Gamma(n + 3 + \alpot)}{\Gamma(n) \, \Gamma(n + 3)} \,
  u_{4, n-1, -1}   \ +  \ 
2 (k_{1} + k_{2})  \, 
\frac{\Gamma(n + 1 + \alpot) \, \Gamma(n + 3 + \alpot)}{\Gamma(n + 1) \, \Gamma(n + 3)} \,
  u_{2, n, -1}
  \right)  \nonumber \\
  & \hspace{3.5in} \ = \ f_{2, n, -1}      \label{erq11}  \\
R6 \, \{ (n , l) \}_{n \ge 1 , l \ge 3}: \quad &   \nonumber  \\
 & \hspace{-1.75in} 2^{\alpha - 2} \, \left( (k_{1} - k_{2}) \, 
\frac{\Gamma(n + \alpot) \, \Gamma(n + 1 + \alpot + l)}{\Gamma(n) \, \Gamma(n + 1 + l)} \,
  u_{l+2, n-1, -1}   \ +  \ 
2 (k_{1} + k_{2})  \, 
\frac{\Gamma(n + 1 + \alpot) \, \Gamma(n + 1 + \alpot + l )}{\Gamma(n + 1) \, \Gamma(n + 1 + l)} \,
  u_{l, n, -1}   \right.    \nonumber  \\
 & \hspace{0.1in} \left. + \  (k_{1} - k_{2}) \, 
\frac{\Gamma(n + 2 + \alpot) \, \Gamma(n + 1 + \alpot + l)}{\Gamma(n + 2) \, \Gamma(n + 1 + l)} \,
  u_{l-2, n+1, -1}
  \right) \ = \ f_{l, n, -1}     \label{erq12}  
\end{align}

\subsection{Determining the coefficients $u_{l, n, \mu}$}
\label{ssec_Detu}
In this section we discuss the solution of $\{ u_{l, n, \mu} \}$ satisfying \eqref{erq1}-\eqref{erq12}. 
We focus our attention on the coefficients $\{ u_{l, n, 1} \}$, i.e. equations  \eqref{erq1}-\eqref{erq6}.
The coefficients $\{ u_{l, n, -1} \}$ are determined in an analogous manner.

Firstly we introduce a rescaling of the equations. Let,
\be
      d_{l, n} \ = \ \frac{\Gamma(n + 1 + \alpot)}{\Gamma(n + 1)} \, u_{l, n, 1} \, , \ \ \mbox{ and } \ \ 
     \wtilde{f}_{l, n} \ = \ 2^{-(\alpha - 2)} \,  \frac{\Gamma(n + 1 + l)}{\Gamma(n + 1 + \alpot + l)} \, f_{l, n, 1} \, .
\label{erq20}
\ee
Then, after simplifying, equations  \eqref{erq1}-\eqref{erq6} become:
\begin{align}
 R1 \, (n = 0 , l = 0):  \quad &  (k_{1} + k_{2}) \,   d_{0, 0}  \ = \ \wtilde{f}_{0, 0}      \label{erq21}  \\
R2 \, (n = 0 , l = 1): \quad &  \left( 2 k_{1} \, + \, (k_{1} + k_{2}) \right) \,  d_{1, 0}  \ = \ \wtilde{f}_{1, 0}      \label{erq22}  \\
R3 \, \{ (n = 0 , l ) \}_{l \ge 2}: \quad & 2 (k_{1} + k_{2}) \,   d_{l, 0}   \ +  \ 
  (k_{1} - k_{2}) \,  d_{l-2, 1}  \ = \ \wtilde{f}_{l, 0}      \label{erq23}  \\
R4 \, \{ (n , l = 0) \}_{n \ge 1}: \quad &     (k_{1} - k_{2}) \,   d_{2, n-1}   \ +  \ 
  (k_{1} + k_{2}) \,  d_{0, n} \ = \ \wtilde{f}_{0, n}      \label{erq24}  \\
R5 \, \{ (n , l = 1) \}_{n \ge 1}: \quad &  (k_{1} - k_{2}) \,  d_{3, n-1}   \ +  \ 
(2 k_{1} \, + \,   (k_{1} + k_{2}) ) \,  d_{1, n} \ = \ \wtilde{f}_{1, n}      \label{erq25}  \\
R6 \, \{ (n , l) \}_{n \ge 1 , l \ge 2}: \quad &  (k_{1} - k_{2}) \,  d_{l+2, n-1}   \ +  \ 
2 (k_{1} + k_{2})  \,   d_{l, n}  \ + \  (k_{1} - k_{2}) \,  d_{l-2, n+1} \ = \ \wtilde{f}_{l, n}    \, .  \label{erq26}  
\end{align}
 
\begin{figure}[!ht]
\begin{center}
 \includegraphics[height=2.0in]{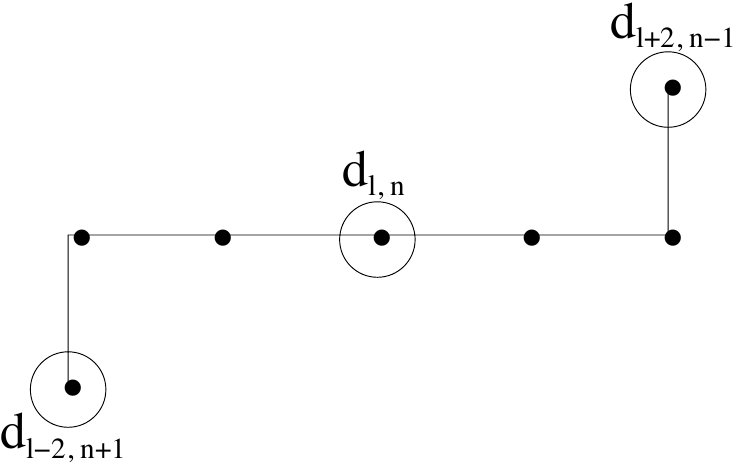}
   \caption{Stencil illustrating the coupling of the unknowns $d_{l , n}$.}
   \label{pic2}
\end{center}
\end{figure}

\underline{The coupling stencil.} \\
Consider the ``stencil'' given in Figure \ref{pic2}. This stencil represents the coupling of the unknowns,
$d_{l, n}$, in equations \eqref{erq21}-\eqref{erq26}. (For the coupling stencil, when $l, \, n < 0$ the $d_{l, n}$
coefficient is ignored.)

Equations \eqref{erq21} and \eqref{erq22} represent two $1\times1$ linear system of equations for $d_{0, 0}$
and $d_{0, 1}$.

From $R3$, $l = 2$ and $n = 0$ gives the equation
\begin{equation}
 2 (k_{1} + k_{2}) \,   d_{2, 0}   \ +  \ 
  (k_{1} - k_{2}) \,  d_{0, 1}  \ = \ \wtilde{f}_{2, 0}     \, .
   \label{erq33}
\end{equation}

Following the coupling of the unknowns from \eqref{erq33} (cf. coefficient stencil), the equation corresponding to 
$l = 0, \, n = 1$ (in $R4$) is
\begin{equation}
 (k_{1} - k_{2})  \,   d_{2, 0}   \ +  \ 
  (k_{1} + k_{2}) \,  d_{0, 1}  \ = \ \wtilde{f}_{0, 1}     \, .
   \label{erq34}
\end{equation}

Returning to $R3$. For $l = 3$ and $n = 0$ we have
\begin{equation}
 2 (k_{1} + k_{2}) \,   d_{3, 0}   \ +  \ 
  (k_{1} - k_{2}) \,  d_{1, 1}  \ = \ \wtilde{f}_{3, 0}     \, .
   \label{erq35}
\end{equation}

Again, following the coupling of the unknowns from \eqref{erq35} (cf. coefficient stencil), the equation corresponding to 
$l = 1, \, n = 1$ (in $R5$) is
\begin{equation}
 (k_{1} - k_{2})  \,   d_{3, 0}   \ +  \ 
(2 k_{1} \, + \,   (k_{1} + k_{2}) \,  d_{1, 1}  \ = \ \wtilde{f}_{1, 1}     \, .
   \label{erq36}
\end{equation}

Equations \eqref{erq33}-\eqref{erq36} are two $2\times2$ linear systems of equations for $d_{2, 0}$, $d_{0, 1}$, and 
$d_{3, 0}$, $d_{1, 1}$, respectively. Continuing in this manner, i.e., in $R3$ corresponding to $(l, n) = (2 m , \, 0)$ 
(and $(l, n) = (2m + 1 , \, 0)$ ), $m \ge 0$, and following the coupling of the equations we terminate in $R4$
corresponding to $(l, n) = (0, \, m)$ 
(in $R5$ with $(l, n) = (1, \, m)$ ) yielding an $(m + 1)\times(m + 1)$ decoupled linear system of equations. Hence, when 
the equations are appropriately assembled, we obtain decoupled linear system of equations of sizes:
1, 1, 2, 2, 3, 3, $\ldots$, $m$, $m$, $\ldots$.

\begin{figure}[!ht]
\begin{minipage}{.46\linewidth}
\begin{center}
 \includegraphics[height=2.25in]{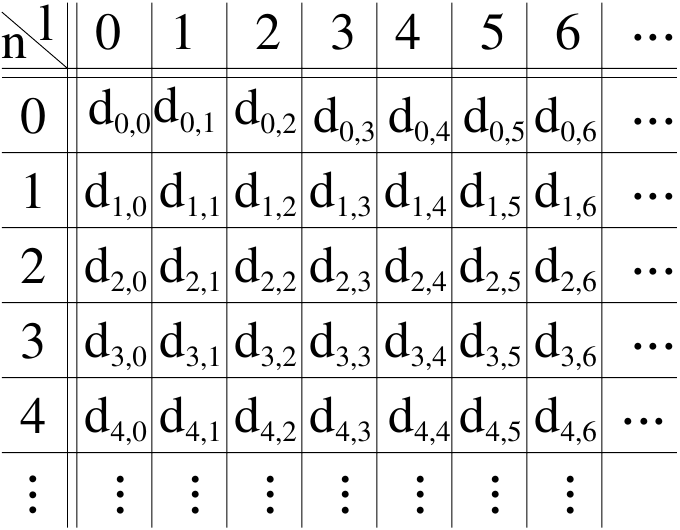}
\end{center}
\end{minipage} \hfill
\begin{minipage}{.46\linewidth}
 
\begin{center}
 \includegraphics[height=2.25in]{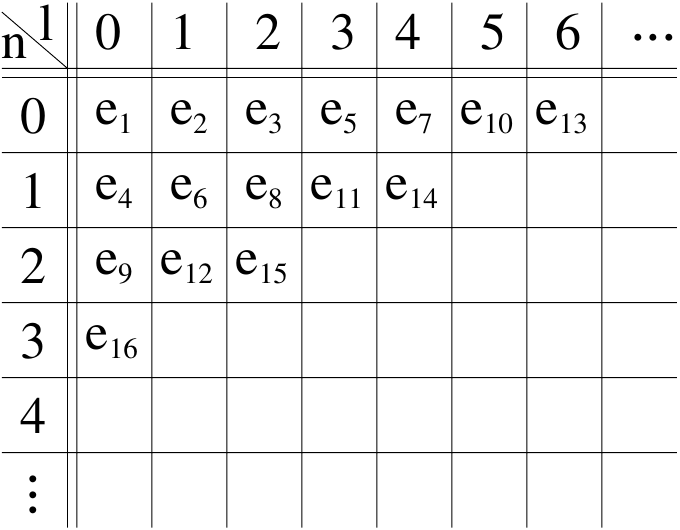}
\end{center}
\end{minipage} 
   \caption{Renumbering of the unknowns $d_{l , n} \rightarrow e_{j}$.}
  \label{renum}
\end{figure}

To further illustrate this point, we rename/renumber the unknowns $d_{l, n} \rightarrow e_{j}$
(and rhs: $\wtilde{f}_{l, n} \rightarrow b_{j}$) in the order described above, see Figure \ref{renum}. 

Algebraically, the rename/renumbering $d_{l, n} \rightarrow e_{j}$ is given by
\begin{align*}
\mbox{if } l \mbox{ is even:} \ & \ (l, n) \longrightarrow j \, = \,  \frac{(l \, + \, 2 n)^{2} \, + \, 2 ( l \, + \, 4 n \, + \, 2)}{4}    \, , \\
\mbox{if } l \mbox{ is odd:} \ & \ (l, n) \longrightarrow j \, = \,   \frac{(l \, + \, 2 n)^{2} \, + \, 2 ( l \, + \, 4 n \, + \, 2) \, + \, 1}{4}     \, .
\end{align*}

As discussed above, the linear system of equations for $\bfe \, = \, \{ e_{j} \}_{j \ge 1}$, $A \bfe \, = \, \bfb$, has the following block
diagonal structure.
\be
A \, = \, \left[ \begin{array}{cccccccc}
k_{1} + k_{2} &  &   &  &  &  &   &    \\
   &   \hspace{-0.5in}   2 k_{1}  +  (k_{1} + k_{2}) &  &   &  &  &  &     \\
   &       A_{3, 4, 2} &  &   &  &  &   &   \\
   &       &  \hspace{-0.3in} A_{3, 5, 2} &  &   &  &    &     \\
   &       &     &  \hspace{-0.3in}   \ddots &  &   &         &   \\
   &       &     &     &  \hspace{-0.3in}  A_{3, 4, m}  &  &       &         \\
   &       &     &     &     &   \hspace{-0.3in}  A_{3, 5,m}  &     &     \\
   &       &     &     &    &     &      \hspace{-0.3in}   \ddots   &   \\
   &       &     &     &    &     &     &       \ddots   
 \end{array}
\right]
\label{Astr}
\ee

The matrices $A_{3, 4, m}$ and $A_{3, 5, m}$ are $m\times m$, symmetric, tridiagonal matrices given by
\be
A_{3, 4, m} \, = \, \left[ \begin{array}{ccccccc}
2(k_{1} + k_{2}) &  (k_{1} - k_{2})&   &  &  &  &    \\
 (k_{1} - k_{2})  &     2(k_{1} + k_{2}) &  (k_{1} - k_{2}) &   &  &  &      \\
   &      &  \hspace{-1.0in} \ddots &  &     &  &        \\
   &       &  & \hspace{-1.0in} \ddots &     &   &           \\
   &       &     &  \hspace{-1.5in} (k_{1} - k_{2})  &   \hspace{-0.5in}   2(k_{1} + k_{2}) &    (k_{1} - k_{2})   &              \\
   &       &     &   &     \hspace{-0.5in}   (k_{1} - k_{2})  &         (k_{1} + k_{2})             \end{array}
\right]
\label{Asub1}
\ee
and
\be
A_{3, 5, m} \, = \, \left[ \begin{array}{ccccccc}
2(k_{1} + k_{2}) &  (k_{1} - k_{2})&   &  &  &  &    \\
 (k_{1} - k_{2})  &     2(k_{1} + k_{2}) &  (k_{1} - k_{2}) &   &  &  &      \\
   &      &  \hspace{-1.0in} \ddots &  &     &  &        \\
   &       &  & \hspace{-1.0in} \ddots &     &   &           \\
   &       &     &  \hspace{-1.5in} (k_{1} - k_{2})  &   \hspace{-0.5in}   2(k_{1} + k_{2}) &    (k_{1} - k_{2})   &              \\
   &       &     &   &     \hspace{-0.5in}   (k_{1} - k_{2})  &         2k_{1} + (k_{1} + k_{2})             \end{array}
\right]
\label{Asub2}
\ee

The matrix $A_{3, 4, m}$ represents the system of equations corresponding to unknowns 
$\bfe_{3,4,m} = \{ e_{j} \} $, $j \, = \, m(m-1) + 1, \ldots, m^{2}$, and
$A_{3, 5, m}$ represents the system of equations corresponding to unknowns $\bfe_{3,5,m} = \{ e_{j} \} $, $j \, = \, m^{2} + 1, \ldots, m(m+1)$.

As $A_{3, 4, m}$ and $A_{3, 5, m}$ are real, symmetric matrices, and hence have real eigenvalues (singular values), a
simple application of Gerschgorin's theorem establishes that there exists constants $c_{min}$ and $c_{max}$ such that
the minimum and maximum eigenvalues of the matrices satisfy $0 < c_{min} \le \lambda_{min}, \, \lambda_{max} \le c_{max} < \infty$.
Hence $A_{3, 4, m}$ and $A_{3, 5, m}$ are uniformly invertible with
\[
  \| \bfe_{3,4,m} \|_{2} \, \le \, \| A_{3, 4, m}^{-1} \|_{2} \, \| \bfb_{3,4,m} \|_{2} \, \le \, c_{min}^{-1}  \, \| \bfb_{3,4,m} \|_{2} \, , 
\]
and, likewise,  $ \| \bfe_{3,5,m} \|_{2} \, \le \,  c_{min}^{-1}  \, \| \bfb_{3,5,m} \|_{2}$.

We summarize the above in the following theorem.
\begin{theorem} \label{thmexst1}
For $k_{1}, k_{2} \in \real^{+}$, with $K(x) \, = \, \left[ \begin{array}{cc}
k_{1}  &   0 \\  0  &  k_{2}  \end{array} \right]$, and 
$f(x) \in L^{2}_{\alpot}(\Omega)$ there exists a unique solution
$\wtilde{u}(x) \ = \ \omega^{\alpot} \, u(x)$,  with $u(x)$ given by \eqref{plk3}, satisfying
\begin{align}
\mcL(\wtilde{u})(x) \, := \, - \Grad \cdot (- \Delta)^{\frac{\alpha - 2}{2}} \, K(x) \, \Grad \wtilde{u}(x) &= \ f(x),  \ \ x \in \Omega\, ,  \label{erq40} \\
\mbox{subject to   }  \wtilde{u} \ = \ 0 &\mbox{ on } \partial \Omega \, ,    \label{erq41} 
\end{align}
where the equal sign in \eqref{erq40} refers to the equality of the coefficients of the basis functions \linebreak[4]
$\{  V_{l , \mu}(x) \, P_{n}^{(\alpot , l)}(2 r^{2} - 1) \}_{l, n, \mu}$.
\end{theorem} 
\mbox{  } \hfill \qed

In the next section we further investigate the relationship between the coefficients $u_{l, n, \mu}$ and $f_{l, n, \mu}$, i.e.,
how the regularity of the solution $\wtilde{u}$ depends on the regularity of the rhs $f$.

We conclude this section with the following corollary that verifies equation \eqref{eqform1} for $\wtilde{g}(x) \ = \ \omega^{\alpot} \, g(x)$, 
with $g(x) \in L^{2}_{\alpot}(\Omega)$.
\begin{corollary} \label{cor4rep}
For $\wtilde{g}(x) \ = \ \omega^{\alpot} \, g(x)$, with $g(x) \in L^{2}_{\alpot}(\Omega)$, we have for $k \in \real$
\be
    k \left(- \Delta \right)^{\alpot} \wtilde{g}(x) \ = \ - \Grad \cdot \left(- \Delta \right)^{\frac{\alpha - 2}{2}} k \bfI \,  \Grad \, \wtilde{g}(x) \, .
\label{eqform1s}
\ee 
\end{corollary}
\textbf{Proof}: Without loss of generality, we assume that $k = 1$. We verify \eqref{eqform1s} by showing that it 
holds on a basis for $L^{2}_{\alpot}(\Omega)$. Specifically, by showing that for $l, \, n \ge 0, \ \mu \in \{ 1 , -1\}$
\[
  \left(- \Delta \right)^{\alpot} \omega^{\alpot} \, V_{l, \mu}(x) \, P_{n}^{(\alpot, l)}(2 r^{2} \, - \, 1) \ = \ 
 - \Grad \cdot \left(- \Delta \right)^{\frac{\alpha - 2}{2}} k \bfI \,  \Grad \omega^{\alpot} \, V_{l, \mu}(x) \, P_{n}^{(\alpot, l)}(2 r^{2} \, - \, 1) \, .
\]
From \cite[Theorem 3]{dyd171},
\begin{align*}
& \left(- \Delta \right)^{\alpot} \omega^{\alpot} \, V_{l, \mu}(x) \, P_{n}^{(\alpot, l)}(2 r^{2} \, - \, 1) 
\ = \ \lambda_{l, n} \,   V_{l, \mu}(x) \, P_{n}^{(\alpot, l)}(2 r^{2} \, - \, 1)  \, ,   \\
& \hspace{0.6in} \mbox{where } \ \  \lambda_{l, n}  \ = \ 2^{\alpha} \, 
\frac{ \Gamma(n + 1 + \alpot) \, \Gamma(n + 1 + \alpot + l) }{ \Gamma(n + 1) \, \Gamma(n + 1 + l) } \, .
\end{align*}

For $l \ge 1, \, n \ge 0, \, \mu \in \{1 , -1\}$
\[
 - \Grad \cdot \left(- \Delta \right)^{\frac{\alpha - 2}{2}} \, \bfI \,  \Grad \omega^{\alpot} \, V_{l, \mu}(x) \, P_{n}^{(\alpot, l)}(2 r^{2} \, - \, 1) 
\ = \ c_{l, n, \mu} \, V_{l, \mu}(x) \, P_{n}^{(\alpot, l)}(2 r^{2} \, - \, 1) \, ,
\]
where $c_{l, n, \mu}$ is the coefficient of $u_{l, n, \mu}$ (for $k_{1} = k_{2} = 1$) in equations \eqref{erq2}, \eqref{erq3},
\eqref{erq5} - \eqref{erq12}, i.e.,
\[
  c_{l, n, \mu}  \ = \ 2^{\alpha - 2} \, 2^{2} \, 
\frac{ \Gamma(n + 1 + \alpot) \, \Gamma(n + 1 + \alpot + l) }{ \Gamma(n + 1) \, \Gamma(n + 1 + l) } 
\ = \ \lambda_{l, n}  \, .
\]

For $l = 0, \, n \ge 0, \, \mu = 1$
\[
 - \Grad \cdot \left(- \Delta \right)^{\frac{\alpha - 2}{2}} \, \bfI \,  \Grad \omega^{\alpot} \, V_{l, \mu}(x) \, P_{n}^{(\alpot, l)}(2 r^{2} \, - \, 1) 
\ = \ 2 \, c_{0, n, 1} \, V_{0, 1}(x) \, P_{n}^{(\alpot, l)}(2 r^{2} \, - \, 1) \, ,
\]
where $c_{0, n, 1}$ is the coefficient of $u_{0, n, 1}$ (for $k_{1} = k_{2} = 1$) in equations \eqref{erq1} and  \eqref{erq4},
 i.e.,
\[
  2 \, c_{0, n, 1}  \ = \ 2 \,  2^{\alpha - 2} \, 2 \, 
\frac{ \Gamma(n + 1 + \alpot) \, \Gamma(n + 1 + \alpot) }{ \Gamma(n + 1) \, \Gamma(n + 1) } 
\ = \ \lambda_{0, n}  \, .
\]
\mbox{  } \hfill \qed

\setcounter{equation}{0}
\setcounter{figure}{0}
\setcounter{table}{0}
\setcounter{theorem}{0}
\setcounter{lemma}{0}
\setcounter{corollary}{0}
\setcounter{definition}{0}
\section{Regularity analysis of the solution}
  \label{sec_Reg}
In Section \ref{sec_ExUn} we established that given $f(x) \in L^{2}_{\alpot}(\Omega)$ 
there exists a 
unique $\wtilde{u}(x) \ = \ \omega^{\alpot} \, u(x)$, with $u(x)$ given by \eqref{plk3}, satisfying
\eqref{erq40},\eqref{erq41}. In this section we analyze the influence of the regularity of $f(x)$ on the
regularity of $u(x)$.

Following \cite{hao211} we introduce the weighted function space $\bfB_{\beta}^{s_{1}, s_{2}}(\Omega)$
for any $s_{1}, \, s_{2} > 0$ as
\[
\bfB_{\beta}^{s_{1}, s_{2}}(\Omega) \ := \ \left\{ v \, | \, v \in L_{\beta}^{2}(\Omega) \ \mbox{ and } \ 
| v |_{\bfB_{\beta}^{s_{1}, s_{2}}(\Omega)} < \infty \right\} \, , 
\]
where the semi-norm $| \cdot |_{\bfB_{\beta}^{s_{1}, s_{2}}(\Omega)}$ is defined by
\[
 | v |^{2}_{\bfB_{\beta}^{s_{1}, s_{2}}(\Omega)} \ = \ 
 \sum_{l, n, \mu} \left( l^{2 s_{1}} \, + \,  n^{2 s_{2}} \right) \, v_{l, n, \mu}^{2} \, h_{l, n}^{2} \, .
\]
The norm in this space is defined as  $\| v \|_{\bfB_{\beta}^{s_{1}, s_{2}}(\Omega)}$, 
where 
\[
 \| v \|^{2}_{\bfB_{\beta}^{s_{1}, s_{2}}(\Omega)} \ = \ 
 \sum_{l, n, \mu} \left( 1 \, + \, l^{2 s_{1}} \, + \,  n^{2 s_{2}} \right) \, v_{l, n, \mu}^{2} \, h_{l, n}^{2} \, .
\]

In determining the coefficients of $u(x)$, $u_{l, n, \mu}$, in Section \ref{sec_ExUn}
the following (invertible) transformations were made:
\be
\left. \begin{array}{l}
      e_{j} \ \longleftarrow \ d_{l, n} \  \longleftarrow \ \quad \frac{\Gamma(n + 1 + \alpot)}{\Gamma(n + 1)} u_{l, n, \mu}  \vspace{1em} \\
      b_{j} \ \longleftarrow \ \wtilde{f}_{l, n} \  \longleftarrow \ 2^{-(\alpha - 2)} \, \frac{\Gamma(n + 1 + l)}{\Gamma(n + 1 + \alpot + l)} f_{l, n, \mu} 
  \end{array}  
  \right\} \ \mbox{ for } \mu = 1 \ \mbox{ and } \mu = 2 \, .
\label{erqq5}
\ee

Similar to the discussion in Section \ref{ssec_Detu}, we focus our attention on $\{ u_{l, n, 1} \}$ and $\{ f_{l, n, 1} \}$, with the
reliance of $\{ u_{l, n, -1} \}$ on $\{ f_{l, n, -1} \}$ following in an analogous manner. To this end, let

\begin{align}
 u_{1}(x) \, = \, \sum_{l \geq 1 , n \geq 0} u_{l, n, 1} \, V_{l , 1}(x) \, P_{n}^{(\alpot , l)}(2 r^{2} - 1) 
&+ \ \sum_{n\geq 0} \frac{u_{0,n,1}}{2} \, V_{0 , 1}(x) \, P_{n}^{(\alpot , 0)}(2 r^{2} - 1) \, , \ \   \label{erqq6} \\
 f_{1}(x) \, = \, \sum_{l \geq 0 , n \geq 0} f_{l, n, 1} \, V_{l , 1}(x) \, P_{n}^{(\alpot , l)}(2 r^{2} - 1)   
 & \ \mbox{ and } \
d(x) \, = \,  \sum_{l \geq 0 , n \geq 0} d_{l, n} \, V_{l , 1}(x) \, P_{n}^{(\alpot , l)}(2 r^{2} - 1) \, , \label{erqq7}
\end{align}
where $\{ d_{l, n} \}$ is determined by \eqref{erq20}-\eqref{erq26}.

We first relate the regularity of $d(x)$ to that of $f_{1}(x)$ and then relate the regularity of $u_{1}(x)$ to $d(x)$ 
(and subsequently to $f_{1}(x)$).

From Section \ref{ssec_Detu} we have that the unknowns $\{ e_{j} \}_{j \ge 1}$ are determined by solving the
sequence of linear systems of equations
\be
A_{3, 4, m} \bfe_{3, 4, m} \, = \, \bfb_{3, 4, m} \, , \ \ \mbox{ and } \ \ A_{3, 5, m} \bfe_{3, 5, m} \, = \, \bfb_{3, 5, m} 
 \, , \ \ \mbox{ for } m = 1, 2, \ldots ,
 \label{erqq10}
\ee
where the $m$x$m$ matrices $A_{3, 4, m}$ and $A_{3, 5, m}$ are symmetric, positive definite with \linebreak[4]
$c_{max}^{-1} \le \| A_{3, 4, m}^{-1} \|_{2} \, , \ \| A_{3, 5, m}^{-1} \|_{2} \le c_{min}^{-1} $. Hence it immediately follows that
\[
\sum_{j = 1}^{\infty} e_{j}^{2} \ \lesssim \sum_{j = 1}^{\infty} b_{j}^{2} \ \ \Longleftrightarrow \ \ 
\sum_{l = 1}^{\infty} \sum_{n = 0}^{\infty} d_{l, n}^{2} \ \lesssim \sum_{l = 1}^{\infty} \sum_{n = 1}^{\infty} \wtilde{f}_{l, n}^{2}
\ = \  \sum_{l = 1}^{\infty} \sum_{n = 0}^{\infty} 
\left( 2^{-(\alpha - 2)} \, \frac{\Gamma(n + 1 + l)}{\Gamma(n + 1 + \alpot + l)} f_{l, n, 1} \right)^{2} \, .
\]

However, to determine which $\bfB_{\alpot}^{s_{1}, s_{2}}(\Omega)$ space $d(x)$ lies in we need to be able to bound
\be
\sum_{l = 1}^{\infty} \sum_{n = 0}^{\infty} (1 \, + \, l^{2 s_{1}} \, + \, n^{2 s_{2}} ) \, d_{l, n}^{2} \, h_{l, n}^{2} \, . 
\label{erqq11}
\ee
Returning to  
\eqref{erqq10}, for $W$ a diagonal matrix, with $w_{i i} > 0$, note that
\[
   W \, A_{3, k, m}  \, W^{-1} \, W  \bfe_{3, k, m} \, = \, W \bfb_{3, k, m} \, , \ \ k \in \{4, 5\} \, .
\]
With $\mcA \, = \,  W \, A_{3, k, m}  \, W^{-1}$, note that $\mcA$ is symmetric, positive definite, and has the same
eigenvalues as $A_{3, k, m}$. Thus, $c_{max}^{-1} \le \| \mcA^{-1} \|_{2} \le c_{min}^{-1}$, and
$ \| W  \bfe_{3, k, m} \|_{2}^{2} \, \lesssim \, \| W  \bfb_{3, k, m} \|_{2}^{2}$. This enables us to introduce the
weights needed in \eqref{erqq11}.

\begin{lemma}  \label{lmareg4d}
For $f_{1}(x)$ and $d(x)$ defined in \eqref{erqq6} and \eqref{erqq7}, if $f_{1}(x) \in \bfB_{\alpot}^{s_{1}, s_{2}}(\Omega)$ then
$d(x) \in \bfB_{\alpot}^{s_{1} + \alpot \, , \, s_{2} + \alpot}(\Omega)$.
\end{lemma}
\textbf{Proof}: From above, we have that
\begin{align*}
\sum_{l = 1}^{\infty} \sum_{n = 0}^{\infty} (1 \, + \, l^{2 t_{1}} \, + \, n^{2 t_{2}} ) \, d_{l, n}^{2} \, h_{l, n}^{2}   
& \lesssim \sum_{l = 1}^{\infty} \sum_{n = 0}^{\infty} (1 \, + \, l^{2 t_{1}} \, + \, n^{2 t_{2}} ) \, \wtilde{f}_{l, n}^{2} \, h_{l, n}^{2}  \\
& = \sum_{l = 1}^{\infty} \sum_{n = 0}^{\infty} (1 \, + \, l^{2 t_{1}} \, + \, n^{2 t_{2}} ) \, 
\left( 2^{-(\alpha - 2)} \, \frac{\Gamma(n + 1 + l)}{\Gamma(n + 1 + \alpot + l)} f_{l, n, 1} \right)^{2}  \, h_{l, n}^{2}  \\
\mbox{(using Stirling's formula)} \ \ \ 
&  \lesssim \sum_{l = 1}^{\infty} \sum_{n = 0}^{\infty} (1 \, + \, l^{2 t_{1}} \, + \, n^{2 t_{2}} ) \, 
\left( (n + 1 + l)^{-\alpot} f_{l, n, 1} \right)^{2}  \, h_{l, n}^{2}      \\
&  \lesssim \sum_{l = 1}^{\infty} \sum_{n = 0}^{\infty} (1 \, + \, l^{2 t_{1}} \, + \, n^{2 t_{2}} ) \, 
 (1 + l + n)^{-\alpha} f_{l, n, 1}^{2}  \, h_{l, n}^{2}      \\
&  \lesssim \sum_{l = 1}^{\infty} \sum_{n = 0}^{\infty} (1 \, + \, l^{2 (t_{1} - \alpot)} \, + \, n^{2 (t_{2} - \alpot)} ) \, 
 f_{l, n, 1}^{2}  \, h_{l, n}^{2}      \, .
\end{align*}
Thus, if $f_{1}(x) \in \bfB_{\alpot}^{s_{1}, s_{2}}(\Omega)$ then
$d(x) \in \bfB_{\alpot}^{s_{1} + \alpot \, , \, s_{2} + \alpot}(\Omega)$.  \\
\mbox{ } \hfill \qed

The connection between the regularity of $u_{1}(x)$ and $d(x)$ is given in the next lemma.
\begin{lemma} \label{lmareg4u1}
For $u_{1}(x)$ and $d(x)$ defined in \eqref{erqq6} and \eqref{erqq7}, if $d(x) \in \bfB_{\alpot}^{s_{1}, s_{2}}(\Omega)$ then
$u_{1}(x) \in \bfB_{\alpot}^{s_{1} \, , \, s_{2} + \alpot}(\Omega)$.
\end{lemma}
\textbf{Proof}: Using \eqref{erqq5}
\begin{align*}
\| u_{1} \|^{2}_{\bfB_{\alpot}^{t_{1}, t_{2}}(\Omega)} & \lesssim \
\sum_{l = 1}^{\infty} \sum_{n = 0}^{\infty} (1 \, + \, l^{2 t_{1}} \, + \, n^{2 t_{2}} ) \, u_{l, n, 1}^{2} \, h_{l, n}^{2}     \\
& \lesssim \sum_{l = 1}^{\infty} \sum_{n = 0}^{\infty} (1 \, + \, l^{2 t_{1}} \, + \, n^{2 t_{2}} ) \, 
\left( \frac{\Gamma(n + 1)}{\Gamma(n + 1 + \alpot)} \, d_{l, n} \right)^{2}  \, h_{l, n}^{2}  \\
\mbox{(using Stirling's formula)} \ \ \ 
&  \lesssim \sum_{l = 1}^{\infty} \sum_{n = 0}^{\infty} (1 \, + \, l^{2 t_{1}} \, + \, n^{2 t_{2}} ) \, 
\left( (n + 1)^{-\alpot} d_{l, n} \right)^{2}  \, h_{l, n}^{2}   \\
& = \sum_{l = 1}^{\infty} \sum_{n = 0}^{\infty} (1 \, + \, l^{2 t_{1}} \, + \, n^{2 t_{2}} ) \, 
(1 + n)^{-\alpha} \,  d_{l, n}^{2}  \, h_{l, n}^{2}     \\
& \lesssim \sum_{l = 1}^{\infty} \sum_{n = 0}^{\infty} (1 \, + \, l^{2 t_{1}} \, + \, n^{2 (t_{2} - \alpot)} ) \, 
 \,  d_{l, n}^{2}  \, h_{l, n}^{2}  \, .
\end{align*}
Thus, if $d(x) \in \bfB_{\alpot}^{s_{1}, s_{2}}(\Omega)$ then
$u_{1}(x) \in \bfB_{\alpot}^{s_{1} \, , \, s_{2} + \alpot}(\Omega)$. \\ 
\mbox{ } \hfill \qed

We are now in the position to state the regularity result for the solution $\tilde{u}(x)$ satisfying \eqref{erq40}, \eqref{erq41}.
\begin{theorem} \label{thmregu1}
For $k_{1}, k_{2} \in \real^{+}$, with $K(x) \, = \, \left[ \begin{array}{cc}
k_{1}  &   0 \\  0  &  k_{2}  \end{array} \right]$, and 
$f(x)  \in \bfB_{\alpot}^{s_{1}, s_{2}}(\Omega)$  there exists a unique solution
$\wtilde{u}(x) \ = \ \omega^{\alpot} \, u(x)$ satisfying
\eqref{erq40}, \eqref{erq41}, with  $u(x) \in \bfB_{\alpot}^{s_{1} + \alpot \, , \, s_{2} + \alpha}(\Omega)$.
\end{theorem} 
\textbf{Proof}: The result follows from Lemmas  \ref{lmareg4d} and  \ref{lmareg4u1} applied to $\{ f_{l, n, 1} \}$, $\{ u_{l, n, 1} \}$, and
$\{ f_{l, n, 2} \}$, $\{ u_{l, n, 2} \}$.
\mbox{  } \hfill \qed

\textbf{Remark}: For the special case of $k_{1} = k_{2}$ Theorem \ref{thmregu1} is consistent with Theorem 3.2 in \cite{hao211}.

\setcounter{equation}{0}
\setcounter{figure}{0}
\setcounter{table}{0}
\setcounter{theorem}{0}
\setcounter{lemma}{0}
\setcounter{corollary}{0}
\setcounter{definition}{0}
\section{Concluding remarks}
  \label{sec_Conc}
In this article we have introduced a generalization of the fractional Laplacian, \eqref{defgenLap}, and showed that for 
$K(x)$ a constant, symmetric positive definite matrix the
associated fractional Poisson equation (in $\real^{2}$) is well posed. Additionally we have related the regularity of the
solution $u(x)$ to the regularity of the RHS function $f(x)$. 
We are currently working to extend these results to a general positive definite matrix $K(x)$, and to extend the results
to $\Omega \subset \real^{3}$. \\

\textbf{Acknowledgements}\\
The third author was partially supported by the National Science Foundation under Grant No. DMS-2012291.


\appendix
 \setcounter{equation}{0}
\setcounter{figure}{0}
\setcounter{table}{0}
\setcounter{theorem}{0}
\setcounter{lemma}{0}
\setcounter{corollary}{0}
\setcounter{definition}{0}
\section{Recurrence formulas for Jacobi polynomials}
 \label{sec_recJ}
In this section we give a number of recurrence formulas for Jacobi polynomials used in the analysis
(see \cite[pg. 782]{abr641}).
\begin{align}
(n \, + \, \alpot \, + \, \frac{\beta}{2} \, + \, 1) \, (1 - x) \, P_{n}^{(\alpha + 1 \, , \beta)}(x) 
&= \ (n \, + \, \alpha \, + \, 1)  \, P_{n}^{(\alpha , \beta)}(x) \ - \   (n \, +  \, 1)  \, P_{n+1}^{(\alpha  , \beta)}(x)   \label{RRJP15} \\
(n \, + \, \alpot \, + \, \frac{\beta}{2} \, + \, 1) \, (1 + x) \, P_{n}^{(\alpha  , \, \beta + 1)}(x) 
&= \ (n \, + \, \beta \, + \, 1)  \, P_{n}^{(\alpha , \beta)}(x) \ + \   (n \, +  \, 1)  \, P_{n+1}^{(\alpha  , \beta)}(x)   \label{RRJP16}  \\
%
(2 n \, + \, \alpha \, + \, \beta) \, P_{n}^{(\alpha -1 \,   ,  \beta)}(x) 
&= \ (n \, + \, \alpha \, + \, \beta)  \, P_{n}^{(\alpha , \beta)}(x) \ - \   (n \, +  \, \beta)  \, P_{n-1}^{(\alpha  , \beta)}(x)   \label{RRJP18} \\
%
(2 n \, + \, \alpha \, + \, \beta) \, P_{n}^{(\alpha  , \, \beta - 1)}(x) 
&= \ (n \, + \, \alpha \, + \, \beta)  \, P_{n}^{(\alpha , \beta)}(x) \ + \   (n \, +  \, \alpha)  \, P_{n-1}^{(\alpha  , \beta)}(x)   \label{RRJP19} 
\end{align}

\begin{lemma} \label{lmaJP1}
For $l \ge 1$ and $n \ge 1$,
\be
l \, P_{n}^{(\gamma , l)}(t) \ + \ \left( \frac{1 + t}{2} \right) (n + \gamma + l + 1) \, P_{n-1}^{(\gamma + 1 \, , \, l + 1)}(t)
\ = \ (n + l) \, P_{n}^{(\gamma + 1 \, , \, l - 1)}(t) \, .
 \label{vjJP1}
\ee
\end{lemma}
\textbf{Proof}:
Using \eqref{RRJP16}
\begin{align*}
& l \, P_{n}^{(\gamma , l)}(t) \ + \  \frac{(n + \gamma + l + 1)}{2} \, (1 + t) \, P_{n-1}^{(\gamma + 1 \, , \, l + 1)}(t)  \\
& \quad  = \  l \, P_{n}^{(\gamma , l)}(t) \ + \ 
  \frac{(n + \gamma + l + 1)}{2 \, (n + \frac{\gamma + 1}{2} + \frac{l}{2})} \,
  \left[ (n + l) \, P_{n-1}^{(\gamma + 1 \, , l)}(t) \ + \ n \, P_{n}^{(\gamma + 1 \, , l)}(t) \right]  \\
& \quad =    \frac{l}{(2 n \, +  \, \gamma + l + 1)} \left[ (n + \gamma + l + 1) \, P_{n}^{(\gamma + 1 \, , l)}(t) \ - \
(n + l) \, \, P_{n-1}^{(\gamma + 1 \, , l)}(t) \right]   \ \mbox{(using \eqref{RRJP18})}  \\
& \quad \quad + \  \frac{(n + \gamma + l + 1)}{(2 n \, +  \, \gamma + l + 1)} \,
  \left[ n \, P_{n}^{(\gamma + 1 \, , l)}(t)  \ + \ (n + l) \, P_{n-1}^{(\gamma + 1 \, , l)}(t) \right]  \\
& \quad = \    \frac{(n + l)}{(2 n \, +  \, \gamma + l + 1)} 
\left[ (n + \gamma + l + 1) \, P_{n}^{(\gamma + 1 \, , l)}(t) \ + \ (n + \gamma + 1) \, P_{n-1}^{(\gamma + 1 \, , l)}(t) \right] \\
& \quad = \    \frac{(n + l)}{(2 n \, +  \, \gamma + l + 1)}  \, (2 n \, +  \, \gamma + l + 1) \, P_{n}^{(\gamma + 1 \, , l - 1)}(t)  
\ \mbox{(using \eqref{RRJP19})}  \\
& \quad  = \ (n + l) \, P_{n}^{(\gamma + 1 \, , l - 1)}(t)  \, .
\end{align*}
\mbox{ } \hfill \qed


\setcounter{equation}{0}
\setcounter{figure}{0}
\setcounter{table}{0}
\setcounter{theorem}{0}
\setcounter{lemma}{0}
\setcounter{corollary}{0}
\setcounter{definition}{0}
\section{Proof of Theorem \ref{genThm3v2}}
 \label{sec_pThm3}
In this section we give a proof of Theorem \ref{genThm3v2}. For brevity, we use the notation used in \cite{dyd171}.

Below the function $G^{m \, n}_{p \, q} \left( \begin{array}{c}
                          \bfa \\  \bfb  \end{array} \, \Big{\vert} \, t \right)$ denotes the Meijer $G$-function, and 
$_{p}\bfF_{q}(\cdot)$ the regularized hypergeometric function.                          

The following five conditions are used in some of the statements involving the Meijer $G$-function.
\begin{align*}
\mbox{Condition S:} & \quad 1 \, - \, \overline{a} \ > \ - \underline{b} \, ,   \\
\mbox{Condition A:} & \quad  p \, + \, q \ <  \ 2 m \, + \, 2 n \, ,   \\
\mbox{Condition B:} & \quad p \, + \, q \ =  \ 2 m \, + \, 2 n \, ,  \ \ p \, = \, q \, ,  \\
\mbox{Condition C:} & \quad p \, + \, q \ =  \ 2 m \, + \, 2 n \, ,  \ \ p \, < \, q \, ,  \\
\mbox{Condition D:} & \quad p \, + \, q \ =  \ 2 m \, + \, 2 n \, ,  \ \ p \, > \, q \, . 
\end{align*}  

The following result is used in the analysis. It summarizes Theorems 1 and 2 in \cite{dyd171}.
\begin{theorem} \cite[Theorem 1 and 2]{dyd171}  \label{DKKThm1}
Let $V(x)$ be a solid harmonic polynomial of degree $l \ge 0$ and $\delta \ = \ d \, + \, 2 l$.
Assume that $-d < \alpha < 0$ or $\alpha > 0$ and parameters $m$, $n$, $p$, $q$, $\bfa$, and $\bfb$
satisfy Condition A, as well as
\[
    2 (1 \, - \, \overline{a}) \ > \ - \alpha \, + \, l \, , \quad -2 \underline{b} \ < \ d \, + \, l \, .
\]
Define $f(x) \ := \ V(x) \, G^{m \, n}_{p \, q} \left( \begin{array}{c}
                          \bfa \\  \bfb  \end{array} \, \Big{\vert} \, | x |^{2} \right)$. Then
\be
(-\Delta)^{\alpha / 2} \, f(x) \ = \ 2^{\alpha} \, V(x) \, 
G^{m+1 \, n+1}_{p+2 \, q+2} \left( \begin{array}{ccc}
        1 \, - \, \frac{\delta + \alpha}{2} \, , &  \bfa  - \alpot \, , & - \alpot \\ 
        0 \, , & \bfb - \alpot \, , & 1 - \frac{\delta}{2} \end{array} \, \Big{\vert} \, | x |^{2} \right)
\label{DKK28}
\ee 
for all $x \neq 0$. The same statement holds under Conditions B, C, and D provided
\[
  2 \overline{\lambda} \ > \ - \alpha \, + \, l \, , \quad 2 \underline{\lambda} \ < \ d \, + \, l \, ;
\]
if Condition B is satisfied, we additionally require that $\nu > 0$ and either $| x | \neq 1$ or
$\nu \ > \ 1 \, + \, \alpha$. 

For the case $\alpha > 0$, the results extends to $x = 0$ whenever both $f$ and 
 $(-\Delta)^{\alpha / 2} \, f$ are continuous at $0$.            \vspace{-1em}              
\end{theorem}
\mbox{  } \hfill    \qed

\textbf{Proof of Theorem \ref{genThm3v2}}: 
Proceeding in a similar manner to the proof of Theorem 3 in \cite{dyd171},
let $V(x) \, = \, V_{l , \mu}(x)$, and using (40a) in  \cite{dyd171} and equation 8.4.49.22 in \cite{pru901}
\begin{align}
f(x) &= \ \frac{\Gamma(n + 1 + \alpot - s)}{\Gamma(n + 1)} \, V(x) \, (1 \, - \, | x |^{2})^{\alpot - s}_{+} \,
_{2}\bfF_{1}\left( \begin{array}{c}
                         -n \ , \ \frac{\delta + \alpha}{2} + n - s \\ 
                                 1 \, + \, \alpot \, - \, s
                             \end{array} \Big{\vert} \, 1 - |x|^{2} \right)   \nonumber \\
&= \   \frac{\Gamma(n + 1 + \alpot - s)}{\Gamma(n + 1)} \, V(x) 
G_{2 \, 2}^{2 \, 0}    \left( \begin{array}{cc}
                         1 + \alpot - s + n \, , &1 - \frac{\delta}{2} - n  \\ 
                              0 \, , &1 - \frac{\delta}{2}
                             \end{array} \Big{\vert} \, |x|^{2} \right)    \nonumber \\
&= \  \frac{\Gamma(n + 1 + \alpot - s)}{\Gamma(n + 1)} \, V(x) 
G_{2 \, 2}^{2 \, 0}    \left( \begin{array}{cc}
                         1 + \alpot - s + n \, , &1 - \frac{\delta}{2} - n  \\ 
                              1 - \frac{\delta}{2} \, , &0
                             \end{array} \Big{\vert} \, |x|^{2} \right)    \ \ ( \mbox{ using the defn. of } G_{2 \, 2}^{2 \, 0} (\cdot)  )  
                             \nonumber \\
&= \  \frac{\Gamma(n + 1 + \alpot - s)}{\Gamma(n + 1)} \, V(x) \, (-1)^{n} \, 
G_{2 \, 2}^{1 \, 1}    \left( \begin{array}{cc}
                         1 - \frac{\delta}{2} - n  \, , &1 + \alpot - s + n \\ 
                              0 \, , &1 - \frac{\delta}{2}
                             \end{array} \Big{\vert} \, |x|^{2} \right)    \ \ \mbox{ using } \cite[(51)]{dyd171} .  \label{bceq2}
\end{align}                             
 
Using Theorem \ref{DKKThm1},
\begin{align}
\lefteqn{ ( - \Delta )^{\frac{\alpha - 2}{2}} V(x) \, G^{1 \, 1}_{2 \, 2}    \left( \begin{array}{cc}
                         1 - \frac{\delta}{2} - n   \, , &  1 + \alpot -s + n \\ 
                              0 \, , & 1 - \frac{\delta}{2}
                             \end{array} \Big{\vert} \, |x|^{2} \right) }  \nonumber \\
&= \ 2^{\alpha - 2} \, V(x) \,  G^{2 \, 2}_{4 \, 4}    \left( \begin{array}{cccc}
 2 - \frac{\delta + \alpha}{2} \, , & 2 - \frac{\delta + \alpha}{2} - n  \, , &  2 - s + n \, , &  1 - \alpot  \\ 
                              0 \, , & 1 - \alpot \, , & 2 - \frac{\delta + \alpha}{2}  \, , & 1 - \frac{\delta}{2}
                             \end{array} \Big{\vert} \, |x|^{2} \right)         \nonumber \\
&= \ 2^{\alpha - 2} \, V(x) \,  G^{2 \, 2}_{4 \, 4}    \left( \begin{array}{cccc}
 2 - \frac{\delta + \alpha}{2} \, , & 2 - \frac{\delta + \alpha}{2} - n  \, , &  2 - s + n \, , &  1 - \alpot  \\ 
                              1 - \alpot \, , & 0 \, , & 2 - \frac{\delta + \alpha}{2}  \, , & 1 - \frac{\delta}{2}
        \end{array} \Big{\vert} \, |x|^{2} \right)     \ \ ( \mbox{ using the defn. of } G^{2 \, 2}_{4 \, 4} (\cdot)  )     \nonumber \\
&= \ 2^{\alpha - 2} \, V(x) \,  G^{1 \, 2}_{3 \, 3}    \left( \begin{array}{ccc}
 2 - \frac{\delta + \alpha}{2} \, , & 2 - \frac{\delta + \alpha}{2} - n  \, , &  2 - s + n   \\ 
                              0 \, , & 2 - \frac{\delta + \alpha}{2}  \, , & 1 - \frac{\delta}{2}
        \end{array} \Big{\vert} \, |x|^{2} \right)     \ \ \mbox{ using } \cite[(22)]{dyd171}    \nonumber \\
&= \ 2^{\alpha - 2} \, V(x) \,  G^{1 \, 2}_{3 \, 3}    \left( \begin{array}{ccc}
 2 - \frac{\delta + \alpha}{2} \, , & 2 - \frac{\delta + \alpha}{2} - n  \, , &  2 - s + n   \\ 
                              0 \, , &  1 - \frac{\delta}{2} \, , &  2 - \frac{\delta + \alpha}{2}
        \end{array} \Big{\vert} \, |x|^{2} \right)     \ \ ( \mbox{ using the defn. of } G^{1 \, 2}_{3 \, 3} (\cdot)  )     \nonumber \\    
&= \ 2^{\alpha - 2} \, V(x) \,  G^{1 \, 1}_{2 \, 2}    \left( \begin{array}{cc}
 2 - \frac{\delta + \alpha}{2} - n  \, , &  2 - s + n   \\ 
                              0 \, , &  1 - \frac{\delta}{2} 
        \end{array} \Big{\vert} \, |x|^{2} \right)     \ \ \mbox{ using } \cite[(21)]{dyd171}     .  \label{bceq3}                                                  
\end{align}

Next, using \cite[(23)]{dyd171}, and \cite[(40b)]{dyd171}
\begin{align}
& G^{1 \, 1}_{2 \, 2}    \left( \begin{array}{cc}
 2 - \frac{\delta + \alpha}{2} - n  \, , &  2 - s + n   \\ 
                              0 \, , &  1 - \frac{\delta}{2} 
        \end{array} \Big{\vert} \, |x|^{2} \right) 
 \ = \ \frac{\Gamma(n \, - \, 1 \, + \, \frac{\delta + \alpha}{2})}{\Gamma(n + 2 - s)} \, 
 _{2}\bfF_{1}\left( \begin{array}{c}
                         n \, - \, 1 \, + \, \frac{\delta + \alpha}{2} \ , \ -1 + s - n  \\ 
                                \frac{\delta}{2}
                             \end{array} \Big{\vert} \,  |x|^{2} \right)   \nonumber \\
&= \ \frac{\Gamma(n \, - \, 1 \, + \, \frac{\delta + \alpha}{2})}{\Gamma(n + 2 - s)} \, 
(-1)^{n + 1 - s}  \, \frac{\Gamma(n  + 2 - s)}{\Gamma(n \, +  1 - s + \, \frac{\delta}{2})} \,                            
     P_{n + 1 - s}^{(\frac{\alpha}{2} - 2 + s \, , \,  \frac{\delta}{2} - 1)}(2 |x|^{2} \, - \, 1) \,   .  \label{bceq4}    
\end{align}

Combining \eqref{bceq2}-\eqref{bceq4} we obtain \eqref{bceq1v}. \\
\mbox{ } \hfill \qed

\setcounter{equation}{0}
\setcounter{figure}{0}
\setcounter{table}{0}
\setcounter{theorem}{0}
\setcounter{lemma}{0}
\setcounter{corollary}{0}
\setcounter{definition}{0}
\section{Proof of Theorem \ref{gradinR2}}
 \label{sec_pThm3p1}
%

\textbf{Proof of Theorem \ref{gradinR2}}:
Recall that in $\real^{2}$, $\Grad \cdot \ = \ \frac{\partial}{\partial r} \cdot \what{\bfr} \ + \ 
\frac{1}{r} \, \frac{\partial}{\partial \varphi} \cdot \what{\bfvarphi} $.

Using \eqref{HLZZe1},
\be
\frac{d}{d r} \left(  (1 - r^{2})^{\alpha/2} \, P_{n}^{(\alpot ,  l)}(2 r^{2} - 1) \right)
\ = \ 
- 2 (n + \alpot) \, r \, (1 - r^{2})^{\alpha/2 - 1} \, P_{n}^{(\alpot-1 \, , \,  l+1)}(2 r^{2} - 1) \, .
\label{bc2eq2}
\ee

Also, 
\be
\frac{\partial}{\partial r} V_{l , \mu}(x) \ = \ l \, r^{-1} \, V_{l , \mu}(x) \, .
\label{bc2eq3}
\ee

Combining \eqref{bc2eq2} and \eqref{bc2eq3} with the definition of $\Grad $ we obtain
\begin{align}
\Grad f(x) 
&= \ 
\left( l \, V_{l , \mu}(x) \, r^{-1} (1 - r^{2})^{\alpha/2} \, P_{n}^{(\alpot ,  l)}(2 r^{2} - 1)   \right.  \nonumber \\
& \quad \quad \quad  \quad \quad \quad  \quad \quad \quad  \left.
\ - \ 2 (n + \alpot) \, V_{l , \mu}(x) \, r \,  (1 - r^{2})^{\alpha/2 - 1} \, P_{n}^{(\alpot-1 \,  ,  \, l+1)}(2 r^{2} - 1)
\right) \what{\bfr}    \nonumber  \\
& \quad + \ 
\left(  \left( \frac{\partial}{\partial \varphi} V_{l , \mu}(x) \right) \, r^{-1} 
   (1 - r^{2})^{\alpha/2} \, P_{n}^{(\alpot ,  l)}(2 r^{2} - 1) \right) \what{\bfvarphi} \, .   \label{bc2eq1}
\end{align}

Using $\bfi \ = \ \cos(\varphi) \what{\bfr} \, - \, \sin(\varphi) \what{\bfvarphi}$, from \eqref{bc2eq1} we obtain
\begin{align}
& \frac{\partial f}{\partial x} \ = \   \nonumber \\
&\cos(\varphi)  \left( l \, V_{l , \mu}(x) \, r^{-1} (1 - r^{2})^{\alpha/2} \, P_{n}^{(\alpot ,  l)}(2 r^{2} - 1)  
\ - \ 2 (n + \alpot) \, V_{l , \mu}(x) \, r \,  (1 - r^{2})^{\alpha/2 - 1} \, P_{n}^{(\alpot-1 \,  ,  \, l+1)}(2 r^{2} - 1) \right) \nonumber \\
& - \ \sin(\varphi) \left(  \left( \frac{\partial}{\partial \varphi} V_{l , \mu}(x) \right) \, r^{-1} 
   (1 - r^{2})^{\alpha/2} \, P_{n}^{(\alpot ,  l)}(2 r^{2} - 1) \right) \, .    \label{ders4} 
\end{align}

With a view to combining the $(1 - r^{2})^{\alpha/2} \, P_{n}^{(\alpot ,  l)}(2 r^{2} - 1)$ terms in \eqref{ders4}, \linebreak[4]
let $H_{1} \ = \ \cos(\varphi)  l \, V_{l , \mu}(x) \, r^{-1} \ - \ 
\sin(\varphi) \left(  \frac{\partial}{\partial \varphi} V_{l , \mu}(x) \right) \, r^{-1}$.\\
Suppose $V_{l , \mu}(x) \, = \, r^{l} \cos(l \varphi)$. 
Then $\frac{\partial}{\partial \varphi} V_{l , \mu}(x) \ = \ - r^{l} \,  l \, \sin(l \varphi)$. In this case,  \vspace{-1em}
\be
H_{1} \ = \ l \, r^{l - 1} \left( \cos(l \varphi) \, \cos(\varphi) \ + \ \sin(l \varphi) \, \sin(\varphi) \right)
\ = \ l \, r^{l - 1}  \cos((l -1) \varphi) \ = \ l \, V_{l-1 \,  , \mu}(x) \, .
\label{ders5}
\ee

Alternatively, suppose $V_{l , \mu}(x) \, = \, r^{l} \sin(l \varphi)$. 
Then $\frac{\partial}{\partial \varphi} V_{l , \mu}(x) \ = \  r^{l} \,  l \, \cos(l \varphi)$, and   \vspace{-1em}
\be
H_{1} \ = \ l \, r^{l - 1} \left( \sin(l \varphi) \, \cos(\varphi) \ - \ \cos(l \varphi) \, \sin(\varphi) \right)
\ = \ l \, r^{l - 1}  \sin((l -1) \varphi) \ = \ l \, V_{l-1 \,  , \mu}(x) \, .
\label{ders6}
\ee

The coefficient of $- \ 2 (n + \alpot) \,  (1 - r^{2})^{\alpha/2 - 1} \, P_{n}^{(\alpot-1 \,  ,  \, l+1)}(2 r^{2} - 1)$
in \eqref{ders4} is 
$H_{2} \, := \, \cos(\varphi) \, V_{l , \mu}(x) \, r$. \\
Suppose $V_{l , \mu}(x) \, = \, r^{l} \cos(l \varphi)$.  Then,   \vspace{-1em}
\begin{align}
H_{2} &= \  r^{l + 1} \,  \cos(l \varphi) \, \cos(\varphi) 
\ = \ \frac{1}{2} \, r^{l + 1} \left(  \cos((l -1) \varphi) \ + \ \cos((l +1) \varphi) \right)  \nonumber \\
&= \ \frac{1}{2} \, r^{2} \, r^{l - 1} \,  \cos((l -1) \varphi) \ + \ \frac{1}{2} \, r^{l + 1}  \,  \cos((l +1) \varphi)  \nonumber \\
&= \ \frac{1}{2} \, r^{2} \, V_{l-1 \, , \mu}(x) \ + \ \frac{1}{2} \, V_{l+1 \, , \mu}(x)  \, .
\label{ders7}
\end{align}
Alternatively, suppose $V_{l , \mu}(x) \, = \, r^{l} \sin(l \varphi)$.  Then,   \vspace{-1em}
\begin{align}
H_{2} &= \  r^{l + 1} \,  \sin(l \varphi) \, \cos(\varphi) 
\ = \ \frac{1}{2} \, r^{l + 1} \left(  \sin((l -1) \varphi) \ + \ \sin((l +1) \varphi) \right)  \nonumber \\
&= \ \frac{1}{2} \, r^{2} \, r^{l - 1} \,  \sin((l -1) \varphi) \ + \ \frac{1}{2} \, r^{l + 1}  \,  \sin((l +1) \varphi)  \nonumber \\
&= \ \frac{1}{2} \, r^{2} \, V_{l-1 \, , \mu}(x) \ + \ \frac{1}{2} \, V_{l+1 \, , \mu}(x)  \, .
\label{ders8}
\end{align}

Rewriting \eqref{ders4} using \eqref{ders5},\eqref{ders6} and \eqref{ders7},\eqref{ders8} we have
\begin{align}
& \frac{\partial f}{\partial x} \ = \   \nonumber \\
& l \, V_{l-1 \, , \mu}(x) \, (1 - r^{2})^{\alpha/2} \, P_{n}^{(\alpot ,  l)}(2 r^{2} - 1)  
\ - \  (n + \alpot) \, V_{l-1 \, , \mu}(x) \, r^{2} \,  (1 - r^{2})^{\alpha/2 - 1} \, P_{n}^{(\alpot-1 \,  ,  \, l+1)}(2 r^{2} - 1)  \nonumber \\
& \ - \  (n + \alpot) \, V_{l+1 \, , \mu}(x) \,  (1 - r^{2})^{\alpha/2 - 1} \, P_{n}^{(\alpot-1 \,  ,  \, l+1)}(2 r^{2} - 1) \, .    \label{ders9} 
\end{align}

Let $H_{3} \ = \  l \,  (1 - r^{2}) \, P_{n}^{(\alpot ,  l)}(2 r^{2} - 1) \ - \ 
 (n + \alpot) \, r^{2} \, P_{n}^{(\alpot-1 \,  ,  \, l+1)}(2 r^{2} - 1) $. With the substitution $t \, = \, 2 r^{2} - 1$,
\begin{align}
H_{3} &= \ \frac{l}{2} \, (1 - t) \, P_{n}^{(\alpot ,  l)}(t) \ - \ \frac{(n + \alpot)}{2} \, (1 + t) \, P_{n}^{(\alpot-1 \,  ,  \, l+1)}(t)
 \nonumber \\
&= \ \frac{l}{(2 n \, + \, \alpot \, + \, l \, + \, 1)} \, \left( (n + \alpot) \,  P_{n}^{(\alpot-1 \,  ,  \, l)}(t) \ - \ 
(n + 1) \, P_{n+1}^{(\alpot-1 \,  ,  \, l)}(t)  \right)  \ \ (\mbox{using } \eqref{RRJP15})  \nonumber  \\
& \quad - \, \frac{(n + \alpot)}{(2 n \, + \, \alpot \, + \, l \, + \, 1)} \, \left( (n + l + 1) \,  P_{n}^{(\alpot-1 \,  ,  \, l)}(t) \ + \ 
(n + 1) \, P_{n+1}^{(\alpot-1 \,  ,  \, l)}(t)  \right)  \ \ (\mbox{using } \eqref{RRJP16})  \nonumber  \\
&= \ - \, \frac{(n + \alpot)}{(2 n \, + \, \alpot \, + \, l \, + \, 1)} \,  (n + 1) \,  P_{n}^{(\alpot-1 \,  ,  \, l)}(t) \ - \ 
\frac{(n + 1)}{(2 n \, + \, \alpot \, + \, l \, + \, 1)}
(n + \alpot + l) \, P_{n+1}^{(\alpot-1 \,  ,  \, l)}(t)   \nonumber  \\
&= \ - \, \frac{(n + 1)}{(2 n \, + \, \alpot \, + \, l \, + \, 1)} \left(  
(n + \alpot + l) \, P_{n+1}^{(\alpot-1 \,  ,  \, l)}(t)  \ + \ (n + \alpot) \,  P_{n}^{(\alpot-1 \,  ,  \, l)}(t) \right)  \nonumber \\
&= \ - \, (n + 1) \, P_{n+1}^{(\alpot-1 \,  ,  \, l-1)}(t)   \ \ (\mbox{using } \eqref{RRJP19}) \, .  \label{ders10}
\end{align}

Substituting \eqref{ders10} into \eqref{ders9} and rearranging we obtain \eqref{ders1}.

The expression for $\frac{\partial f}{\partial y}$ is obtained in a similar manner. Using $\bfj \ = \ \sin(\varphi) \what{\bfr} \, + \, \cos(\varphi) \what{\bfvarphi}$, from \eqref{bc2eq1} we obtain
\begin{align}
& \frac{\partial f}{\partial y} \ = \   \nonumber \\
&\sin(\varphi)  \left( l \, V_{l , \mu}(x) \, r^{-1} (1 - r^{2})^{\alpha/2} \, P_{n}^{(\alpot ,  l)}(2 r^{2} - 1)  
\ - \ 2 (n + \alpot) \, V_{l , \mu}(x) \, r \,  (1 - r^{2})^{\alpha/2 - 1} \, P_{n}^{(\alpot-1 \,  ,  \, l+1)}(2 r^{2} - 1) \right) \nonumber \\
& + \ \cos(\varphi) \left(  \left( \frac{\partial}{\partial \varphi} V_{l , \mu}(x) \right) \, r^{-1} 
   (1 - r^{2})^{\alpha/2} \, P_{n}^{(\alpot ,  l)}(2 r^{2} - 1) \right) \, .    \label{ders14} 
\end{align}

With a view to combining the $(1 - r^{2})^{\alpha/2} \, P_{n}^{(\alpot ,  l)}(2 r^{2} - 1)$ terms in \eqref{ders14}, \linebreak[4]
let $H_{1} \ = \ \sin(\varphi)  l \, V_{l , \mu}(x) \, r^{-1} \ + \ 
\cos(\varphi) \left(  \frac{\partial}{\partial \varphi} V_{l , \mu}(x) \right) \, r^{-1}$.\\
Suppose $V_{l , \mu}(x) \, = \, r^{l} \cos(l \varphi)$. 
In this case,  \vspace{-1em}
\be
H_{1} \ = \ l \, r^{l - 1} \left( \cos(l \varphi) \, \sin(\varphi) \ - \ \sin(l \varphi) \, \cos(\varphi) \right)
\ = \ - \,  l \, r^{l - 1}  \sin((l -1) \varphi) \ = \ - \, l \, V_{l-1 \,  , \mu^{*}}(x) \, .
\label{ders15}
\ee

Alternatively, suppose $V_{l , \mu}(x) \, = \, r^{l} \sin(l \varphi)$. 
Then   \vspace{-1em}
\be
H_{1} \ = \ l \, r^{l - 1} \left( \sin(l \varphi) \, \sin(\varphi) \ + \ \cos(l \varphi) \, \cos(\varphi) \right)
\ = \ l \, r^{l - 1}  \cos((l -1) \varphi) \ = \ l \, V_{l-1 \,  , \mu^{*}}(x) \, .
\label{ders16}
\ee

The coefficient of $- \ 2 (n + \alpot) \, (1 - r^{2})^{\alpha/2 - 1} \, P_{n}^{(\alpot-1 \,  ,  \, l+1)}(2 r^{2} - 1)$
in \eqref{ders4} is 
$H_{2} \, := \, \sin(\varphi) \, V_{l , \mu}(x) \, r$. \\
Suppose $V_{l , \mu}(x) \, = \, r^{l} \cos(l \varphi)$.  Then,   \vspace{-1em}
\begin{align}
H_{2} &= \  r^{l + 1} \,  \cos(l \varphi) \, \sin(\varphi) 
\ = \ - \, \frac{1}{2} \, r^{l + 1} \left(  \sin((l -1) \varphi) \ - \ \sin((l +1) \varphi) \right)  \nonumber \\
&= \ - \, \frac{1}{2} \, r^{2} \, r^{l - 1} \,  \sin((l -1) \varphi) \ + \ \frac{1}{2} \, r^{l + 1}  \,  \sin((l +1) \varphi)  \nonumber \\
&= \ - \, \frac{1}{2} \, r^{2} \, V_{l-1 \, , \mu^{*}}(x) \ + \ \frac{1}{2} \, V_{l+1 \, , \mu^{*}}(x)  \, .
\label{ders17}
\end{align}
Alternatively, suppose $V_{l , \mu}(x) \, = \, r^{l} \sin(l \varphi)$.  Then,   \vspace{-1em}
\begin{align}
H_{2} &= \  r^{l + 1} \,  \sin(l \varphi) \, \sin(\varphi) 
\ = \  \frac{1}{2} \, r^{l + 1} \left(  \cos((l -1) \varphi) \ - \ \cos((l +1) \varphi) \right)  \nonumber \\
&= \ \frac{1}{2} \, r^{2} \, r^{l - 1} \,  \cos((l -1) \varphi) \ - \ \frac{1}{2} \, r^{l + 1}  \,  \cos((l +1) \varphi)  \nonumber \\
&= \ \frac{1}{2} \, r^{2} \, V_{l-1 \, , \mu^{*}}(x) \ - \ \frac{1}{2} \, V_{l+1 \, , \mu^{*}}(x)  \, .
\label{ders18}
\end{align}

Rewriting \eqref{ders14} using \eqref{ders15},\eqref{ders16} and \eqref{ders17},\eqref{ders18} we have
\begin{align}
& \frac{\partial f}{\partial y} \ = \   \nonumber \\
& l \, (\mp) V_{l-1 \, , \mu^{*}}(x) \, (1 - r^{2})^{\alpha/2} \, P_{n}^{(\alpot ,  l)}(2 r^{2} - 1)  
\ - \  (n + \alpot) \, (\mp) V_{l-1 \, , \mu^{*}}(x) \, r^{2} \,  (1 - r^{2})^{\alpha/2 - 1} \, P_{n}^{(\alpot-1 \,  ,  \, l+1)}(2 r^{2} - 1)  \nonumber \\
& \ - \  (n + \alpot) \, (\pm) V_{l+1 \, , \mu^{*}}(x) \,  (1 - r^{2})^{\alpha/2 - 1} \, P_{n}^{(\alpot-1 \,  ,  \, l+1)}(2 r^{2} - 1) \, .    \label{ders19} 
\end{align}

Using \eqref{ders10} and rearranging we obtain \eqref{ders2}.  The proofs of (\ref{ders1z0}) and (\ref{ders2z0}) follow in a similar and much simpler manner (as  $V_{0,1}=1$), and are thus omitted.
\mbox{  } \hfill \qed



\begin{thebibliography}{10}

\bibitem{abr641}
M.~Abramowitz and I.A. Stegun.
\newblock {\em Handbook of mathematical functions with formulas, graphs, and
  mathematical tables}, volume~55 of {\em National Bureau of Standards Applied
  Mathematics Series}.
\newblock For sale by the Superintendent of Documents, U.S. Government Printing
  Office, Washington, D.C., 1964.

\bibitem{aco171}
G.~Acosta and J.P. Borthagaray.
\newblock A fractional {L}aplace equation: {R}egularity of solutions and finite
  element approximations.
\newblock {\em SIAM J. Numer. Anal.}, 55(2):472--495, 2017.

\bibitem{aco181}
G.~Acosta, J.P. Borthagaray, O.~Bruno, and M.~Maas.
\newblock Regularity theory and high order numerical methods for the
  (1-d)-fractional {L}aplacian.
\newblock {\em Math. Comp.}, 87:1821--1857, 2018.

\bibitem{caf071}
L.~Caffarelli and L.~Silvestre.
\newblock An extension problem related to the fractional {L}aplacian.
\newblock {\em Comm. Partial Differential Equations}, 32(7-9):1245--1260, 2007.

\bibitem{dun141}
C.F. Dunkl and Y.~Xu.
\newblock {\em Orthogonal polynomials of several variables}, volume 155 of {\em
  Encyclopedia of Mathematics and its Applications}.
\newblock Cambridge University Press, Cambridge, second edition, 2014.

\bibitem{dyd171}
B.~Dyda, A.~Kuznetsov, and M.~Kwa\'{s}nicki.
\newblock Fractional {L}aplace operator and {M}eijer {G}-function.
\newblock {\em Constr. Approx.}, 45(3):427--448, 2017.

\bibitem{erv191}
V.J. Ervin.
\newblock Regularity of the solution to fractional diffusion, advection,
  reaction equations in weighted {S}obolev spaces.
\newblock {\em J. Differential Equations}, 278:294--325, 2021.

\bibitem{hao211}
Z.~Hao, H.~Li, Z.~Zhang, and Z.~Zhang.
\newblock Sharp error estimates of a spectral {G}alerkin method for a
  diffusion-reaction equation with integral fractional {L}aplacian on a disk.
\newblock {\em Math. Comp.}, 90(331):2107--2135, 2021.

\bibitem{hao201}
Z.~Hao and Z.~Zhang.
\newblock Optimal {R}egularity and {E}rror {E}stimates of a {S}pectral
  {G}alerkin {M}ethod for {F}ractional {A}dvection-{D}iffusion-{R}eaction
  {E}quations.
\newblock {\em SIAM J. Numer. Anal.}, 58(1):211--233, 2020.

\bibitem{kwa171}
M.~Kwa\'{s}nicki.
\newblock Ten equivalent definitions of the fractional {L}aplace operator.
\newblock {\em Fract. Calc. Appl. Anal.}, 20(1):7--51, 2017.

\bibitem{li141}
H.~Li and Y.~Xu.
\newblock Spectral approximation on the unit ball.
\newblock {\em SIAM J. Numer. Anal.}, 52(6):2647--2675, 2014.

\bibitem{li211}
Y.~Li.
\newblock On the decomposition of solutions: from fractional diffusion to
  fractional {L}aplacian.
\newblock {\em Fract. Calc. Appl. Anal.}, 24(5):1571--1600, 2021.

\bibitem{lis201}
A~Lischke, G.~Pang, M.~Gulian, F.~Song, C.~Glusa, X.~Zheng, Z.~Mao, W.~Cai,
  M.M. Meerschaert, M.~Ainsworth, and G.Em Karniadakis.
\newblock What is the fractional {L}aplacian? {A} comparative review with new
  results.
\newblock {\em J. Comput. Phys.}, 404:109009, 62, 2020.

\bibitem{mao161}
Z.~Mao, S.~Chen, and J.~Shen.
\newblock Efficient and accurate spectral method using generalized {J}acobi
  functions for solving {R}iesz fractional differential equations.
\newblock {\em Appl. Numer. Math.}, 106:165--181, 2016.

\bibitem{mao181}
Z.~Mao and G.Em Karniadakis.
\newblock A spectral method (of exponential convergence) for singular solutions
  of the diffusion equation with general two-sided fractional derivative.
\newblock {\em SIAM J. Numer. Anal.}, 56(1):24--49, 2018.

\bibitem{pru901}
A.P. Prudnikov, Yu.A. Brychkov, and O.I. Marichev.
\newblock {\em Integrals and series. {V}ol. 3}.
\newblock Gordon and Breach Science Publishers, New York, 1990.
\newblock More special functions, Translated from the Russian by G. G. Gould.

\bibitem{zen141}
F.~Zeng, F.~Liu, C.~Li, K.~Burrage, I.~Turner, and V.~Anh.
\newblock A {C}rank-{N}icolson {ADI} spectral method for a two-dimensional
  {R}iesz space fractional nonlinear reaction-diffusion equation.
\newblock {\em SIAM J. Numer. Anal.}, 52(6):2599--2622, 2014.

\bibitem{zhe211}
X.~Zheng, V.J. Ervin, and H.~Wang.
\newblock Optimal {P}etrov-{G}alerkin spectral approximation method for the
  fractional diffusion, advection, reaction equation on a bounded interval.
\newblock {\em J. Sci. Comput.}, 86(3):Paper No. 29, 22, 2021.

\end{thebibliography}

\end{document}